\documentclass[number,citesort,MSNbibl,dvips]{arxbj}
\usepackage{shere}
\usepackage[left]{vmkcol}
\usepackage{graphicx}

\aid{0}
\volume{18}
\issue{1}
\pubyear{2012}
\firstpage{290}
\lastpage{321}
\doi{10.3150/10-BEJ338}

\makeatletter
\newcolumntype{d}[1]{D{.}{.}{#1}}

\def\a{\alpha}
\def\pa{\mathrm{pa}}

\def\cC{\mathcal{C}}
\def\cG{\mathcal{G}}
\def\cH{\mathcal{H}}
\def\cX{\mathcal{X}}
\def\cU{\mathcal{U}}
\def\cL{\mathcal{L}}
\def\cM{\mathcal{M}}
\def\cK{\mathcal{K}}
\def\mm{\mathfrak{m}}
\def\R{\mathbb{R}}
\def\N{\mathbb{N}}
\def\P{\mathbb{P}}
\def\E{\mathbb{E}}

\newtheorem{theo}{Theorem}[section]
\newtheorem{lem}[theo]{Lemma}
\newtheorem{prop}[theo]{Proposition}
\newtheorem{cor}[theo]{Corollary}
\newproclaim{defn}[theo]{Definition}
\newremark{exmp}[theo]{Example}
\newremark{rem}[theo]{Remark}

\newcommand\indep{\protect\mathpalette{\protect\independenT}{\perp}}
\def\independenT#1#2{\mathrel{\rlap{$#1#2$}\mkern4.1mu{#1#2}}}

\renewcommand{\emptyset}{\varnothing}
\makeatother

\begin{document}
\begin{frontmatter}

\title{Tree cumulants and the geometry of binary tree models}
\runtitle{Tree cumulants and the geometry of binary tree models}

\begin{aug}
\author{\fnms{Piotr} \snm{Zwiernik}\corref{}\thanksref{e1}\ead[label=e1,mark]{piotr.zwiernik@gmail.com}} \and
\author{\fnms{Jim Q.} \snm{Smith}\thanksref{e2}\ead[label=e2,mark]{j.q.smith@warwick.ac.uk}}
\runauthor{P.~Zwiernik and J.Q.~Smith}
\address{Department of Statistics, University of Warwick, Coventry CV4 7AL,
UK.\\
\printead{e1,e2}}
\end{aug}

\received{\smonth{5} \syear{2010}}
\revised{\smonth{10} \syear{2010}}

%
\begin{abstract}
In this paper we investigate undirected discrete graphical tree models
when all the variables in the system are binary,
where leaves represent the observable variables and where all the inner
nodes are unobserved. A~novel approach based on the theory of
partially ordered sets allows us to obtain a convenient parametrization
of this model class. The construction of the proposed coordinate
system mirrors the combinatorial definition of cumulants. A simple
product-like form of the resulting parametrization gives insight
into identifiability issues associated with this model class. In
particular, we provide necessary and sufficient conditions for such
a model to be identified up to the switching of labels of the inner
nodes. When these conditions hold, we give explicit formulas for
the parameters of the model. Whenever the model fails to be identified,
we use the new parametrization to describe the geometry
of the unidentified parameter space. We illustrate these results using
a simple example.
\end{abstract}

%

\end{frontmatter}

\section{Introduction}

Discrete graphical models have become a very popular tool in the
statistical analysis of multivariate problems (see, e.g., \cite
{lauritzen96,spiegelhalter1993bae}). When all the variables in
the system are observed, they exhibit a useful modularity. In
particular, it is possible to estimate all the conditional
probabilities that parametrize such models, maximum likelihood
estimates are simple sample proportions and a conjugate Bayesian
analysis is straightforward. However, if the values of some of the
variables are unobserved, then the resulting model for the observed
variables often becomes very complex, making inference much more difficult.

The complicated structure of models with hidden variables usually leads
to difficulties in establishing the identifiability of their parameters
(see, e.g., \cite{allman2009identifiability}). In this paper, we show
how algebraic and combinatorial techniques can help. We focus on
graphical models where the underlying graph is a tree and all the inner
nodes represent hidden variables. In the computational biology
literature, these models are called the general Markov models (see,
e.g., \cite{semple2003pol}), tree models or tree decomposable
distributions (cf. \cite{pearltarsi86}). Building on results of Chang
\cite{chang1996frm}, in this paper we analyze issues associated with
identifiability of such a tree model when all its variables are binary,
paying particular attention to the geometry of the unidentified space.
In particular, we obtain necessary and sufficient conditions for this
model to be locally identified, which gives a stronger version of
Theorem 4.1 in \cite{chang1996frm}. When these conditions are
satisfied, we also obtain exact formulae for its parameters in terms of
the marginal distribution over the observed variables.

Our strategy is to define a new parametrization of this model class.
The new coordinate system is based on moments rather than conditional
probabilities. This helps us to exploit various invariance properties
of tree models, which, in turn, enables us to express the dependence
structure implied by the tree more elegantly. Furthermore, because the
parametrization is based on well-understood moments, the implied
dependence structure becomes more transparent.

\begin{figure}

\includegraphics{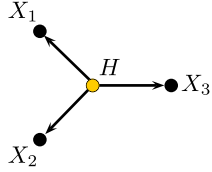}%
\vspace*{-3pt}
\caption{The tripod tree model.}\label{fig:tripod}
\vspace*{-4pt}
\end{figure}

The motivation of this methodology sprung from the study of the tripod
tree model, which is the simplest naive Bayes model. The model is a
graphical model given in Figure~\ref{fig:tripod},
where the black nodes represent three observed variables,
$X_{1},X_{2},X_{3}$, and the white node indicates a hidden variable $H$
that remains \textit{hidden}; that is, its values are never directly observed.
We assume all the variables in the system have values in $\{0,1\}$. For
$\a=(\a_{1},\a_{2},\a_{3})\in\{0,1\}^{3}$ let $p_{\a}=\P
(X_{1}=\a_{1},X_{2}=\a_{2},X_{3}=\a_{3})$. This model would usually
be parametrized using conditional probabilities. In this case we would write
\begin{equation}\label{eq:param-tripod}
p_{\a}=\sum_{i=0}^{1}\theta^{(h)}_{i} \theta^{(1)}_{\a
_{1}|i}\theta^{(2)}_{\a_{2}|i}\theta^{(3)}_{\a_{3}|i},
\end{equation}
where $\theta^{(h)}_{i}=\P(H=i)$ and $\theta^{(j)}_{\a_{j}|i}=\P
(X_{j}=\a_{j}|H=i)$. It can be seen that there are\vspace*{-2pt} seven free
parameters needed to specify $p_{\a}$, namely: $\theta^{(h)}_{1}$
together with $\theta^{(j)}_{1|i}$ for $i=0,1$ and $j=1,2,3$.

However, the definition of this model given in (\ref{eq:param-tripod})
becomes more transparent when expressed in terms of moments. It is easy
to check that there is a one-to-one correspondence between the
probabilities $p_{\a}$ for $\a\in\{0,1\}^{3}$ and the four central
moments $\mu_{ij}=\E(X_{i}-\lambda_{i})(X_{j}-\lambda_{j})$ {for }
$i,j=1,2,3$ and
$\mu_{123}= \E(X_{1}-\lambda_{1})(X_{2}-\lambda_{2})(X_{3}-\lambda
_{3})$ supplemented by the three means $\lambda_{i}=\E X_{i}$ {for }
$i=1,2,3$ (cf. Appendix \ref{sec:change}).

Let $\bar{\mu}_{h}=1-2\theta_{1}^{(h)}$, $\bar{\mu}_{i}=1-2\lambda
_{i}$ and $\eta_{h,i}=\theta_{1|1}^{(i)}-\theta_{1|0}^{(i)}$ for
$i=1,2,3$. We can now\vspace*{-2pt} write an explicit isomorphism\vadjust{\goodbreak} between the
original seven parameters $(\theta^{(h)}_{1},(\theta
^{(i)}_{1|0},\theta^{(i)}_{1|1}))$ and\vspace*{2pt} new parameters $(\bar{\mu
}_{h},(\bar{\mu}_{i}), (\eta_{h,i}))$ for $i=1,2,3$. Thus, in \cite
{settimi1998gbg}, it is shown that in the new coordinate system,
together with the new parameters, the model class is equivalently given by
\begin{eqnarray}\label{eq:star}
\lambda_{i}&=&\tfrac{1}{2}(1-\bar{\mu}_{i})\qquad\mbox{for }
i=1,2,3,\nonumber\\[-2pt]
\mu_{ij}&=&\tfrac{1}{4}(1-\bar{\mu}_{h}^2)\eta_{h,i}\eta_{h,j} \qquad \mbox
{for all } i\neq j\in\{1,2,3\}  \quad \mbox{and}\\[-2pt]
\mu_{123}&=&\tfrac{1}{4}(1-\bar{\mu}_{h}^2)\bar{\mu}_{h} \eta
_{h,1}\eta_{h,2}\eta_{h,3}.\nonumber
\end{eqnarray}

The product-like form of this parametrization enables us to see various
interesting constraints on the observed nodes. For example, by
multiplying formulae for $\mu_{12},\mu_{13}$ and $\mu_{23}$ in (\ref
{eq:star}) together we can see that $\mu_{12}\mu_{13}\mu_{23}\geq0$
must hold. It also allows us to find explicit formulae for the
parameters of the model in terms of the marginal distribution on the
set of observed variables. For example, when $\mu_{12}\mu_{13}\mu
_{23}\neq0$ by substituting (\ref{eq:star}) for all the observed
moments, we see that
\begin{equation}\label{eq:params-tripod}
\bar{\mu}_{h}^{2}=\frac{\mu_{123}^{2}}{\mu_{123}^{2}+4\mu_{12}\mu
_{13}\mu_{23}},\qquad\eta_{h,i}^{2}=\frac{\mu_{123}^{2}+4\mu
_{12}\mu_{13}\mu_{23}}{\mu_{jk}^{2}} \qquad \mbox{for } i=1,2,3.
\end{equation}

Now a similar parametrization is known for general naive Bayesian
models; see the Appendix in \cite{geiger2001sef}. The new
parametrization for this model class was used in \cite{rusakov2006ams}
to approximate a marginal likelihood where the sample size was large,
in \cite{geiger2001sef} to understand the local geometry of the model
class and in \cite{auvray2006sad} to provide the full description of
these models in terms of the defining equations and inequalities.

Naive Bayesian models are a particular example of general Markov
models. The class of tree models is somewhat more complicated than the
naive Bayesian models and needs new tools to examine its geometry. In
this paper, we investigate the moment structures induced by tree models
using the theory of partially ordered sets and M\"{o}bius functions.
Similar methods were used in the combinatorial theory of cumulants (see
\cite{rotacumulants,speed1983cumulants}) for a poset of all
partitions of a finite set. To our knowledge, this paper is the first
to use more general posets in statistical analysis, although a similar
approach can be found in the theory of free probability (see, e.g.,
\cite{speicher94}).


The paper is organized as follows. In Section \ref{sec:models_trees}
we define and analyze the moment structures of the class of models
under consideration. In Section \ref{sec:monomial} we define
tree-cumulants, which form a new coordinate system for this model
class. In Section \ref{sec:parametrization} we reparametrize the model
and show that the induced parametrization on the observed margin has an
elegant product-like form. We apply this reparametrization in Section
\ref{sec:fibers}, analyzing the local geometry of the tree models and
the geometry of the subsets of the parameter space that give the same
set of marginal distributions on the set of observed variables. In
Section \ref{sec:quartet} we illustrate this method using a simple
general Markov model given by a~tree with two hidden nodes.

\vspace*{-3pt}\section{Independence models on trees}\vspace*{-3pt}\label{sec:models_trees}

In this section, we introduce models defined by global Markov
properties on trees.\vadjust{\goodbreak}

\subsection{Preliminaries on trees}\label{sec:intro_trees}

A \textit{graph} $G$ is an ordered pair $(V,E)$ consisting of a
non-empty set $V$ of \textit{nodes} (or \textit{vertices}) and a set
$E$ of \textit{edges}, each of which is an element of $V\times V$. An
edge $(u,v)\in E$ is \textit{directed} if the pair $(u,v)$ is ordered
and we represent the edge by an arrow from $u$ to $v$. If $(u,v)$ is
not an ordered pair, then we say that $(u,v)$ is an \textit{undirected
edge}. Graphs with only (un)directed edges are called (un)directed. If
$e=(u,v)$ is an edge of a graph $G$, then $u$ and $v$ are called
\textit{adjacent} and $e$ is said to be \textit{incident with} $u$
and $v$. If $v\in V$, the \textit{degree} of~$v$ is denoted by $\deg
(v)$, and is the number of edges incident with $v$. A \textit{path} in
a graph~$G$ is a sequence of nodes $(v_1,v_2,\ldots, v_{k})$ such
that, for all $i=1,\ldots,k-1$, $v_i$ and $v_{i+1}$ are adjacent. If,
in addition, $v_1=v_k$, then the path is called a \textit{cycle}. A
graph is \textit{connected} if each pair of nodes in $G$ can be joined
by a path.

A (\textit{directed}) \textit{tree} $T=(V,E)$ is a connected (\textit
{directed}) graph with no cycles. A node of~$T$ of degree one is called
a \textit{leaf}. A node of $T$ that is not a leaf is called an \textit
{inner node}. An edge $e$ of $T$ is \textit{inner} if both nodes
incident with $e$ are inner nodes. A connected subgraph of $T$ is a
\textit{subtree} of~$T$. A \textit{rooted tree}, $T^{r}$, is a
directed tree that has one distinguished node called the \textit
{root}, denoted by the letter $r$, and edges that are directed away
from $r$. Let $T^{r}$ be a rooted tree. For every node $v$ of $T^{r}$
we let $\pa(v)$ denote the set of nodes $u$ such that $(u,v)\in E$. If
$v$ is the root, then $\pa(v)=\emptyset$. Otherwise $\pa(v)$ is a singleton.\looseness=1

For any $W\subseteq V$ we define $T(W)$ as the minimal subtree of $T$
whose set of nodes contains~$W$. We say that $T(W)$ is the subtree of
$T$ \textit{spanned on $W$}. Henceforth, denote the edge set of $T(W)$
by $E(W)$ and its set of nodes by $V(W)$. If $T$ is rooted, then let~$r(W)$ denote the unique node $v$ of $T(W)$ such that $\pa(v)\cap
V(W)$ is the empty set.

Let $T=(V,E)$ be a tree where $e=(u,v)$ denotes one of its edges. Then
\textit{contracting}~$e$ results in another tree, denoted by $T/e$,
with the edge $e$ removed and its incident nodes~$u$ and $v$
identified. Similarly, for any $E'\subseteq E$ we denote the tree
obtained from $T$ by contracting all edges in $E'$ by $T/E'$. If $v\in
V$ such that $\deg v=2,$ then to \textit{suppress $v$} we simply
contract one of the edges incident with $v$. The resulting tree is
denoted by $T/ v$.

\subsection{Models defined by global Markov properties}\label{sec:subs-indep}

In this paper, we always assume that random variables are binary,
taking either value~$0$ or $1$. The vector $Y$ has as its components
all variables in the graphical model, that is, both hidden and observed
variables. Denote the subvector of $Y$ of observed variables by~$X$ and
the subvector of {hidden} variables by~$H$.

Let $T=(V,E)$ be an undirected tree. For any three disjoint subsets
$A,B,C\subseteq V$ we say that $C$ \textit{separates} $A$ and $B$ in
$T$, denoted by $A\perp_T B|C$, if each path from a node in~$A$ to a
node $B$ passes through a node in $C$. For any $A\subseteq V$ let
$Y_{A}$ denote the subvector of $Y=(Y_v)_{v\in V}$ with elements
indexed by $A$, that is, $Y_A=(Y_v)_{v\in A}$. We are interested in
statistical models for $Y$ defined by global Markov properties (GMP) on
$T$. By definition (see, e.g., \cite{lauritzen96}, Section 3.2.1),
these models are specified through the set of conditional independence
statements of the form:
\begin{equation}\label{eq:GMP}
 \{Y_A\indep Y_B|Y_C\dvt \mbox{for all } A,B,C\subset V\mbox{
s.t. } A\perp_T B|C \}.
\end{equation}
Let $\widetilde{\cM}_{T}$ denote the space of probability
distributions of $(X,H)$ satisfying the global Markov properties on
$T$. We now let $\cM_{T}$ denote the space of marginal probability
distribution on $X$ induced from distributions over $(X,H),$ which are
in $\widetilde{\cM}_{T}$.

\subsection{Models for rooted trees}\label{sec:rooted}

We next present the parametric formulation of the models presented in
the previous section. A~Markov process on a rooted tree $T^r$ is a
collection of random variables, \mbox{$\{Y_v\dvt v\in V\},$} such that for each
$\a=(\a_{v})_{v\in V}\in\{0,1\}^{V}$
\begin{equation}\label{eq:p_albar}
p_\a(\theta)=\prod_{v\in V} \theta^{(v)}_{\a_{v}|\a_{\pa(v)}},
\end{equation}
where $\pa(r)$ is the empty set, $\theta=(\theta^{(v)}_{\a_{v}|\a
_{\pa(v)}})$ and
\[
\theta^{(v)}_{\a_{v}|\a_{\pa(v)}}=\P\bigl(Y_v={\a}_v|Y_{\pa(v)}=\a
_{\pa(v)}\bigr).
\]
Since $\theta^{(r)}_{0}+\theta^{(r)}_{1}=1$ and $\theta
^{(v)}_{0|i}+\theta^{(v)}_{1|i}=1$ for all $v\in V\setminus\{r\}$ and
$i=0,1$, the set of parameters consists of exactly $2|E|+1$ free
parameters: we have two parameters, $\theta^{(v)}_{1|0}$, $\theta
^{(v)}_{1|1}$, for each edge $(u,v)\in E$ and one parameter, $\theta
^{(r)}_{1}$, for the root. We denote the parameter space by $\Theta
_{T}=[0,1]^{2|E|+1}$.

Suppose that $T^{r}$ has $n$ leaves representing a binary random
vector, $X=(X_{1},\ldots,X_{n}),$ and let
\begin{equation}\label{eq:prob-simplex}
\Delta_{2^n-1}=\biggl\{p\in\R^{2^n}\dvt \sum_\beta p_\beta=1, p_\beta\geq
0\biggr\}
\end{equation}
with indices $\beta$ ranging over $\{0,1\}^n$ be the probability
simplex of all possible distributions of $X$. Equation (\ref
{eq:p_albar}) induces a polynomial map, $f_T\dvtx\Theta_{T}\rightarrow
\Delta_{2^n-1}$, obtained by marginalization over all the inner nodes
of $T$, giving the marginal mass function $p_{\beta}(\theta)$ as
\begin{equation}\label{eq:p_albar2}
p_{\beta}(\theta)=\sum_\cH\prod_{v\in V} \theta^{(v)}_{\a_{v}|\a
_{\pa(v)}}.
\end{equation}
Here, $\cH$ denotes the set of all $\a\in\{0,1\}^{V}$ such that the
restriction to the leaves of $T$ is equal to $\beta$. The image of
this map is, by definition, the general Markov model on $T^{r}$ (cf.~\cite{semple2003pol}, Section~8.3, \cite{pearltarsi86}).

Standard theory in graphical models tells us that the Markov process on
$T^{r}$ is equal to $\widetilde{\cM}_T$ and, consequently, that the
general Markov on $T^{r}$ model is equal to $\cM_{T}$. Indeed, since
$T^{r}$ is a perfect directed graph (see Section 2.1.3 in \cite
{lauritzen96}), by \cite{lauritzen96}, Theorem 3.28, the Markov
properties are equivalent to the factorization with respect to the
undirected version of $T^{r}$, which is just $T$. Since $T$ is
decomposable, by \cite{lauritzen96}, Proposition 3.19, the
factorization according to $T$ is equivalent to the global Markov
properties on $T$.

In this paper, we often focus on \textit{trivalent trees}, that is,
trees such that every inner node has degree three. This is an important
subclass because, by the well-known lemma below (see, e.g., \cite
{pearltarsi86}, Section 2), the nodes of valency two in a given tree
add nothing to the model class $\cM_{T}$.
\begin{lem}\label{rem:deg2edge}
Let $T$ be a tree. Let $v\in V$ be a node of degree two and let
$T'=T/v$ be the tree obtained from $T$ by suppressing $v$. Then $P\in
{\cM}_T$ if and only if $P\in{\cM}_{T'}$.
\end{lem}

\begin{cor}\label{cor:marg-trip}
Let $T$ be a tree and let $i,j,k$ be any three leaves of $T$. The
marginal model on $(X_{i},X_{j},X_{k})$ induced from $\cM_{T}$ and
denoted by $\cM_{T(ijk)}$ is equivalent to the tripod tree model where
the tripod tree is given in Figure \ref{fig:tripod}.
\end{cor}


In addition, the model corresponding to any tree is a submodel of a
model corresponding to a trivalent tree. To show this, we need the
following definition.
\begin{defn}\label{def:triv-exp}
Let $T$ be any tree. A \textit{trivalent expansion} of $T$, denoted by
${T}^{*}$, is any tree ${T}^{*}=({V}^{*},{E}^{*})$ whose each inner
node has degree at most three and there exists a set of inner nodes
$E'\subseteq{E}^{*}$ such that $T={T}^{*}/E'$.
\end{defn}
\begin{lem}\label{lem:triv-red}
Let $T$ be a tree and ${T}^{*}=({V}^{*},{E}^{*})$ its trivalent
expansion with $E'\subseteq E^{*}$ such that $T=T^{*}/E'$. Then $\cM
_{T}\subseteq\cM_{{T}^{*}}$.
\end{lem}
\begin{pf}
Let $p$ be a point in $\cM_{T}$. Then $p=f_{T}(\theta)$ for some
$\theta\in\Theta_{T}$. Identifying edges of~${T}^{*}$ and $T$ in the
obvious way, we can write ${E}^{*}=E'\cup E$. Define ${\theta}^{*}\in
\Theta_{{T}^{*}}$ as follows. For all $\a_{u},\a_{v}\in\{0,1\}$
\begin{eqnarray}\label{eq:star-constraints}
{\theta^{*}}^{(v)}_{\a_{v}|\a_{u}}&=&\theta^{(v)}_{\a_{v}|\a_{u}}  \qquad
\mbox{for every }(u,v)\in E,\nonumber
\\[-8pt]
\\[-8pt]
{\theta^{*}}^{(v)}_{\a_{v}|\a_{u}}&=&\delta_{\a_{u}\a_{v}}  \qquad  \mbox
{for every }(u,v)\in E',\nonumber
\end{eqnarray}
where $\delta_{ij}$ denotes the Kronecker's delta. It is now simple to
check that $f_{{T}^{*}}({\theta^{*}})=p$. It follows that $p\in\cM_{T^{*}}$.
\end{pf}

For these reasons, we can usually safely restrict our attention to
trivalent trees.

\subsection{Moments and conditional independence}\label
{sec:correlations}\label{sec:corr_model}

Let $X=(X_1,\ldots, X_n)$ be a random vector and for each $\beta
=(\beta_1,\ldots,\beta_n)\in\N^n$ denote $X^\beta=\prod_i
X_i^{\beta_i}$. We shall denote $\E X^\beta$ by $\lambda_\beta$ and
$\E U^{\beta}$ by $\mu_\beta$, where $U_{i}=X_{i}-\E X_{i}$. When
$\beta\in\{0,1\}^{n}$, it is often convenient to use an alternate
notation. Thus, for subsets $I\subseteq[n]:=\{1,2,\ldots,n\},$ we let
$\lambda_{I}=\E(\prod_{i\in I} X_{i})$, $\mu_{I} = \E(\prod_{i\in
I} U_{i})$. Note that $\lambda_{e_{i}}$, where $e_{i}$ is the standard
basis vector in $\R^{n}$, can also be denoted by $\lambda_{i}$ for
$i=1,\ldots,n$.

The model $\cM_{T}$ in the previous section is given in terms of the
probabilities as the image of the map in (\ref{eq:p_albar2}). We find
it convenient to change these coordinates. Let $[n]_{\geq2}$ denote
all subsets of $[n]$ with at least two elements. Denote by $\cC_{n}$
the set of values of all the means $\lambda_{1},\ldots, \lambda_{n}$
together with central moments $\mu_{I}$ such that $I\in[n]_{\geq2}$
for all possible probabilities in $\Delta_{2^{n}-1}$. There exists a
polynomial isomorphism, \mbox{$f_{p\mu}\dvtx\Delta_{2^{n}-1}\rightarrow\cC
_{n}$}, with the inverse denoted by $f_{\mu p}$ (for details see
Appendix \ref{sec:change}). Consequently, we can express any
distribution in the general Markov model in terms of its central
moments and means.\looseness=1

For any two sets $A,B$ let $AB$ denote $A\cup B$. If ${X}_A\indep
{X}_B$, then $\mu_{IJ}=\mu_I\mu_J$ for all non-empty $I\subseteq A$,
$J\subseteq B$. However, when all variables are binary, we also have a
converse result. Thus, if for all non-empty $I\subseteq A$, $J\subseteq
B$ we have that $\mu_{IJ}=\mu_I\mu_J$, then ${X}_A\indep{X}_B$.
Indeed, the independence expressed in terms of moments (see, e.g.,
Feller \cite{feller1971ipt}, page 136) gives
\begin{equation}\label{eq:ind-mom}
{X}_A\indep{X}_B \quad \Longleftrightarrow \quad {\operatorname{Cov}}(f(X_A), g(X_B))=0
\qquad \mbox{for all } f\in L^{2}(\cX_{A}),g\in L^{2}(\cX_{B}).
\end{equation}
Since our variables are binary, all the functions of $X_A$ and $X_B$
are just polynomials with square-free monomials. Equivalently, every
function of $X_{A}$ or $X_{B}$ can be written as a~polynomial with
square-free monomials in $U_{A}$ or $U_{B}$, respectively. For
instance, because~$X_{1}, X_{2}\in\{0,1\},$
\[
X_{1}^{10}X_{2}^{3}=X_{1}X_{2}=(U_{1}+\lambda_{1})(U_{2}+\lambda
_{2})=U_{1}U_{2}+\lambda_{2}U_{1}+\lambda_{1}U_{2}+\lambda
_{1}\lambda_{2}.
\]
Since the covariance is a bilinear form, Settimi and Smith \cite
{settimi2000gma} concluded that the independence can be checked only on
these monomials and (\ref{eq:ind-mom}) can be rewritten as
\begin{equation}\label{eq:ind-mom2}
{X}_A\indep{X}_B \quad \Longleftrightarrow \quad {\operatorname{Cov}}(U_{A}^{\a},
U_{B}^{\beta})=0 \qquad \mbox{for all } \a\in\{0,1\}^{|A|}, \beta\in\{
0,1\}^{\beta}.
\end{equation}
However, $\operatorname{{Cov}}(U_A^{\a}, U_B^\beta)=0$ holds for each non-zero
$\a\in\{0,1\}^{|A|}$ and $\beta\in\{0,1\}^{|B|}$ if and only if
$\mu_{IJ}=\mu_{I}\mu_{J}$ for each $I\subseteq A$, $J\subseteq B$.




We can generalize the result above. For a random variable $H_{a}$ let
$\lambda_{a}=\E H_{a}$ and $U_{a}=H_{a}-\lambda_{a}$. For each
$I\subseteq[n]$ let $U_{I}=\prod_{i\in I} U_{i}$ and
\begin{equation}\label{eq:eta-dajesz}
\eta_{a,I}=\E(U_{I}U_{a})/\operatorname{{Var}}(H_{a}).
\end{equation}
Note that under this notation $\operatorname{{Var}}(H_{a})=\lambda_{a}(1-\lambda_{a})$.

\begin{prop}\label{prop:indep}
Let $H_{a}$ be a non-degenerate random variable. With the notation
above, we have $X_A\indep X_B|H_{a}$ if and only if for all non-empty
$I\subseteq A$, $J\subseteq B$
\begin{eqnarray}\label{eq:gen_rhoxy|h}
\mu_{IJ}&=&\mu_I\mu_J+\lambda_{a}(1-\lambda_{a})\eta_{a,I}\eta
_{a,J},\nonumber
\\[-8pt]
\\[-8pt]
\eta_{a,IJ}&=&\mu_I\eta_{a,J}+\eta_{a,I}\mu_J+(1-2\lambda_{a})\eta
_{a,I}\eta_{a,J}.
\nonumber
\end{eqnarray}
\end{prop}
\begin{pf}
The definition of independence given in (\ref{eq:ind-mom2}) induces a
condition for $X_A\indep X_B|H_{a}$. Thus, for each $I\subseteq A$,
$J\subseteq B$ we have
\begin{equation}\label{eq:cond_covs}
\operatorname{{Cov}}(U_{I},U_{J}|H_{a}=0)=\operatorname{{Cov}}(U_{I},U_{J}|H_{a}=1)=0,
\end{equation}
so, in particular,
\begin{eqnarray} \label{eq:cond_covs2}
\lambda_{a}
\operatorname{{Cov}}(U_I,U_J|H_{a}=1)+(1-\lambda_{a})\operatorname{{Cov}}(U_I,U_J|H_{a}=0)&=&0,\nonumber
\\[-9pt]
\\[-9pt]
\operatorname{Cov}(U_I,U_J|H_{a}=0)-\operatorname{Cov}(U_I,U_J|H_{a}=1)&=&0.
\nonumber
\end{eqnarray}
Moreover, for any $I\subseteq[n],$ one has $\E(U_{I}|H_{a})=\mu
_{I}+\eta_{a,I}U_{a}$, and hence
\begin{equation}\label{eq:cond_covs3}
\operatorname{Cov}(U_I,U_J|H_{a})=\mu_{IJ}-\mu_{I}\mu_{J}+(\eta_{a,IJ}-\eta
_{a,I}\mu_{J}-\mu_{I}\eta_{a,J})U_{a}-\eta_{a,I}\eta_{a,J}U_{a}^{2}.
\end{equation}
Equation (\ref{eq:gen_rhoxy|h}) now follows from substituting (\ref
{eq:cond_covs3}) into (\ref{eq:cond_covs2}).
\end{pf}

\section{Tree posets and tree cumulants}\label{sec:monomial}

In this section, we use the theory of partially ordered sets to propose
a further change of coordinates. In the new coordinate system it is
possible to parametrize the marginal model $\cM_{T}$ in a product form
(see Proposition \ref{prop:monomial}) in contrast to the complicated
polynomial mapping given in~(\ref{eq:p_albar2}).


\subsection{The poset of edge partitions}

Let $T=(V,E)$ be a tree with $n$ leaves. We identify the set of leaves
of $T$ with the set~$[n]$. For any $e\in E$ we let $T\setminus e$
denote the \textit{forest} obtained from $T$ by removing $e$, that is,
the subgraph of $T$ given as a collection of disjoint trees with the
set of nodes given by $V$ and the set of edges given by $E\setminus e$.
Similarly, for any $E'\subseteq E$, we let $T\setminus E'$ denote the
forest obtained by removing all the edges in $E'$. An \textit{edge
split} is a partition of the set of leaves,~$[n]$, of $T$ into two
non-empty sets induced by removing an edge $e$ from $E$ and
restricting~%
$[n]$ to the connected components of $T\setminus e$. By an \textit
{edge partition}, we mean any partition $B_1|B_{2}|\cdots|B_k$ of the
set of leaves induced by considering connected components of
$T\setminus E'$ for some $E'\subseteq E$. Call each subset $B_i$ in
this partition a \textit{block}.

Henceforth let $\Pi_{T}$ denote the poset of all edge partitions of
the set of leaves induced by edges of $T$. The ordering is induced from
the ordering of the poset of all partitions of the set of leaves (see
\cite{stanley2006enumerative}, Example 3.1.1.d). Thus, for two
partitions, $\pi=B_{1}|\cdots|B_{k}$ and $\nu=C_{1}|\cdots|C_{l}$,
we write $\pi\leq\nu$ if every block of $\pi$ is contained in one
of the blocks of $\nu$. To make this more explicit, define the
following equivalence relation on the subsets of $E$. For
$E_{1},E_{2}\subseteq E$ we say $E_{1}\sim E_{2}$ if and only if
removing $E_{1}$ induces the same partition of the set of leaves $[n]$
as removing $E_{2}$. For example, in Figure \ref{fig:tripod} the
partition, $1|2|3,$ can be obtained either by removing any two edges or
by removing all them. However, the only way to obtain the partition,
$12|3,$ is by removing the edge incident with the third leaf.

Let $\overline{E}_{\pi}$ denote the element of the equivalence class
of subsets of $E$ inducing the partition $\pi$, which is maximal with
respect to inclusion. Suppose that $\pi\in\Pi_{T}$ is obtained by
removing edges in the subset of the set of edges $E_{\pi}$ and $\nu
\in\Pi_T$ is obtained by removing edges in $E_{\nu}$. Write $\pi
\leq\nu$ if and only if $\overline{E}_{\pi}\supseteq\overline
{E}_{\nu}$ and call $\pi$ a \textit{subpartition} of $\nu$.

An \textit{interval}, $[\pi,\nu]$, for $\pi$ and $\nu$ in $\Pi
_T$, is the set of all elements $\delta$ such that $\pi\leq\delta
\leq\nu$. The poset $\Pi_T$ forms a lattice (cf.\vadjust{\goodbreak} \cite
{stanley2006enumerative}, Section 3.3). To show this, we define $\pi
\vee\nu\in\Pi_{T}$ ($\pi\wedge\nu\in\Pi_{T}$) as an element in
$\Pi_T$ obtained by removing $\overline{E}_\pi\cap\overline{E}_\nu
$ ($\overline{E}_\pi\cup\overline{E}_\nu$). We have $\pi\vee\nu
\geq\pi$, $\pi\vee\nu\geq\nu$ ($\pi\wedge\nu\leq\pi$, $\pi
\wedge\nu\leq\nu$) and, if there exists another $\delta\in\Pi_T$
such that $\delta\geq\pi$, $\delta\geq\nu$ ($\delta\leq\pi$,
$\delta\leq\nu$), then $\delta\geq\pi\vee\nu$ ($\delta\leq\pi
\wedge\nu$). The element $\pi\vee\nu$ ($\pi\wedge\nu$) is
called the \textit{join} (the \textit{meet}) of $\pi$ and $\nu$.
The poset $\Pi_{T}$ has a unique minimal element, $1|2|\cdots|n,$
induced by removing all edges in $E$ and the maximal one with no edges
removed, which is equal to a single block, $[n]$. The maximal and
minimal element of a lattice will be denoted by $\hat{1}$ and $\hat
{0}$, respectively.


The number of elements in these posets is typically large. However, the
key concepts can be presented using a simpler poset. Let $\widetilde
{\Pi}_{T}$ denote a subposet of $\Pi_{T}$ containing partitions
obtained by removing only inner edges and consider, for example, the
two different trivalent trees $T$ and $T'$, both with six leaves, given below
\begin{center}
\vspace*{5pt}
\includegraphics{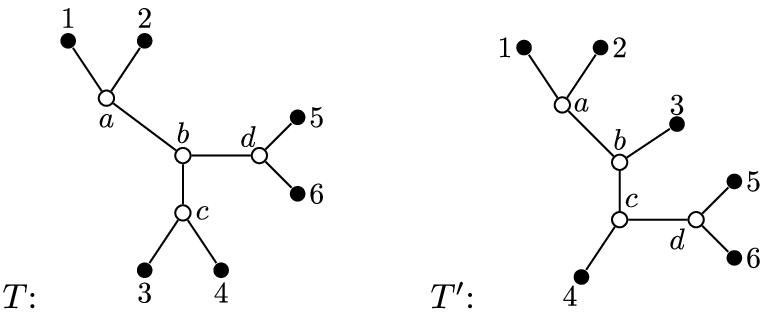}
\vspace*{5pt}
\end{center}
Their associated posets, $\widetilde{\Pi}_T$ and $\widetilde{\Pi
}_{T'}$, are, respectively,
\begin{center}
\vspace*{5pt}
\includegraphics{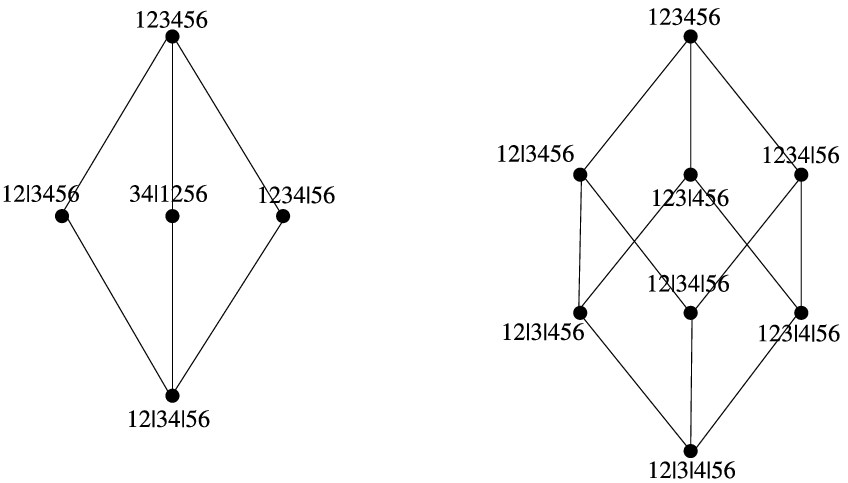}
\vspace*{5pt}
\end{center}

So, for example, $12|34|56$ is an edge partition in $\widetilde{\Pi
}_T$ and is a subpartition of any other edge partition $\nu\in
\widetilde{\Pi}_T$. It can be obtained by removing either any two
inner edges from $(a,b)$, $(b,c)$ and $(b,d),$ or all of them. Since,
for $\pi= 12|34|56,$ there are no subpartitions of~$\pi,$ it follows
that $\pi$ is the minimal element of $\widetilde{\Pi}_T$. In
$\widetilde{\Pi}_{T'}$, there is only one way to obtain this
partition. Namely, by removing $(a,b)$ and $(c,d)$. However, note that
this partition is not minimal in $\widetilde{\Pi}_{T'}$ because
$12|3|4|56< \pi$.

For any poset $\Pi$ a \textit{M\"{o}bius function} $\mm_\Pi\dvtx\Pi
\times\Pi\rightarrow\R$ is defined by $\mm_\Pi(\pi,\pi)=1$ for
every $\pi\in\Pi$,
$\mm_\Pi(\pi,\nu)=-\sum_{\pi\leq\delta<\nu} \mm_\Pi(\pi
,\delta)$ for $\pi<\nu$ in $\Pi$ and is zero otherwise (cf.~\cite
{stanley2006enumerative}, Section~3.7). Recall that\vadjust{\goodbreak} for any $W\subset
V$, $T(W)$ denotes the subtree of $T$ spanned on~$W$ (see Section~\ref
{sec:intro_trees}). We denote $\mm_{\Pi_{T(W)}}:=\mm_W$ and $\mm
_{\Pi_T}:=\mm$, and let $\hat{0}_W$ and $\hat{1}_W$ denote the
minimal and the maximal element of $\Pi_{T(W)}$, respectively. For any
partition $\pi\in\Pi_T$ the interval $[\hat{0},\pi]$ has a natural
structure of a product of posets for blocks of $\pi$, namely~$\prod
_{B\in\pi} \Pi_{T(B)}$, where the product is over all blocks $B$ of
$\pi$. By Proposition 3.8.2 in \cite{stanley2006enumerative}, the M\"{o}bius
function on the product of posets $\prod_{B\in\pi} \Pi
_{T(B)}$ can be written as a~product of M\"{o}bius functions for each
of the posets $ \Pi_{T(B)}$. Thus, for $\nu\leq\pi$ in $\Pi_{T}$
\begin{equation}\label{eq:mob_prod}
\mm(\nu,\pi)=\prod_{B\in\pi} \mm_B(\nu_B,\hat{1}_B),
\end{equation}
where $\nu_B\in\Pi_{T(B)}$ is the restriction of $\nu\in\Pi_T$ to
the block containing only elements from $B\subset[n]$ (it is well
defined since $\nu\leq\pi$) and $\pi_B=\hat{1}_B$ for each $B$.

In the next section, we will use the M\"{o}bius function of the poset
of tree partitions to derive a useful change of coordinates on $\cM_{T}$.

\subsection{An induced change of coordinates}\label{sec:kappas}
Assume that each inner node of $T$ has degree at most three and
consider a map, $f_{\mu\kappa}\dvtx \R^{n}\times\R^{2^{n}}\rightarrow
\R^{n}\times\R^{2^{n}}$, where the coordinates in the domain are
denoted by $\lambda_{1},\ldots, \lambda_{n}$ and $\mu_{I}$ for
$I\subseteq[n]$ and the coordinates in the image are denoted by
$\lambda_{1},\ldots, \lambda_{n}$ and $\kappa_{I}$ for $I\subseteq
[n]$. The map is defined as the identity on the first $n$ coordinates
corresponding to the means and
\begin{equation}\label{eq:kappa-in-rho}
\kappa_{I}=\sum_{\pi\in\Pi_{{T}(I)}}\mm_I(\pi,\hat{1}_{I})
\prod_{B\in\pi}\mu_B \qquad\mbox{for all } I\subseteq[n].
\end{equation}
It is easy to prove that the Jacobian of $f_{\mu\kappa}$ is equal to
$1,$ so, in particular, this is constant. To see this, order the
variables in such a way that the first $n$ coordinates both in $\cK
_{T}$ and~$\cC_{n}$ are $\lambda_{1},\ldots,\lambda_{n}$ and let
$\kappa_I$ precede $\kappa_J$ ($\mu_I$ precede $\mu_J$) as long as
$I\subset J$. The Jacobian matrix of $f_{\mu\kappa}$ is then lower
triangular with each of its diagonal entries equal to~$1$. It follows
that the modulus of its determinant is always $1$.

The map, $f_{\mu\kappa}$, is a regular polynomial map with a regular
polynomial inverse $f_{\kappa\mu}$. Therefore, it gives a change of
coordinates from the central moments with means to a~coordinate system
given by $\lambda_{1},\ldots,\lambda_{n}$ and $\kappa_{I}$ for
$I\subseteq[n]$. Its inverse map is given by\looseness=1
\begin{equation}\label{eq:muinkappa}
\mu_{I}=\sum_{\pi\in\Pi_{{T}(I)}}\prod_{B\in\pi} \kappa_B\qquad
\mbox{for all } I\in[n]_{\geq2}.
\end{equation}\looseness=0
To show (\ref{eq:muinkappa}), define two functions on $\Pi_{{T}(I)}$:
$\a(\pi)=\prod_{B\in\pi} \mu_B$ and $\beta(\pi)=\prod_{B\in
\pi} \kappa_B$. For each $\pi\in\Pi_{{T}(I)}$, by (\ref{eq:kappa-in-rho}),
\begin{eqnarray*}
\beta(\pi) & = & \prod_{B\in\pi} \kappa_B=\prod_{B\in\pi}
\biggl(\sum_{\nu_B\in\Pi_{{T}(B)}}\mm_B(\nu_B,\hat{1}_B)\prod_{C\in
\nu_B}\mu_C  \biggr) \\
& = & \sum_{\nu\leq\pi} \prod_{B\in\pi} \mm_B(\nu_B,\hat
{1}_B) \prod_{C\in\nu}\mu_C,
\end{eqnarray*}
where $\nu$ is an element of $\Pi_{{T}(I)}$ such that its restriction
to each of the blocks $B\in\pi$ is equal to $\nu_{B}$. By the
product formula in (\ref{eq:mob_prod}), we have $\prod_{B\in\pi}
\mm_B(\nu_B,\hat{1}_B)=\mm_{I}(\nu,\pi)$. Therefore, $\beta(\pi
)=\sum_{\nu\leq\pi} \mm_I(\nu,\pi)\a(\nu)$ for all $\pi\in
\Pi_{{T}(I)}$. Equation (\ref{eq:muinkappa}) now follows on applying
the M\"{o}bius inversion formula in Proposition 3.7.1 in \cite
{stanley2006enumerative}.

Denote $\cK_{T}=f_{\mu\kappa}(\cC_{n})$. Since $\cK_{T}$ is
contained in a subset of $\R^{n}\times\R^{2^{n}}$ given by $\kappa
_{\emptyset}=\kappa_{1}=\cdots=\kappa_{n}=0,$ a system of
coordinates on $\cK_{T}$ is given by $\lambda_{i}$ for $i=1,\ldots,
n$ and $\kappa_{I}$ for $I\in[n]_{\geq2}$. This system of
coordinates is called \textit{tree cumulants}. The name is justified
by~(\ref{eq:kappa-in-rho}) because one of the definitions of classical
cumulants is the following. Let $\Pi(I)$ denote the set of all
partitions of $I=\{i_{1},\ldots, i_{k}\}\in[n]_{\geq2}$ (see \cite
{stanley2006enumerative}, Example 3.1.1.d). Then, for all $k>1$
\begin{equation}\label{eq:cumulant}
\operatorname{Cum}(X_{i_{1}},\ldots, X_{i_{k}})=\sum_{\pi\in\Pi(I)} \mm
_{\Pi(I)}(\pi,\hat{1}_{I})\prod_{B\in\pi}\mu_{B},\vspace*{-2pt}
\end{equation}
where the product is over all blocks of $\pi$. Moreover, for every
$\pi\in\Pi(I)$
\[
\mm_{\Pi(I)}(\pi,\hat{1}_{I})=(-1)^{|\pi|-1}(|\pi|-1)!,\vspace*{-2pt}
\]
where $|\pi|$ denotes the number of blocks in $\pi$. Note that the
usual definition of cumulants uses non-central moments instead of
central moments in (\ref{eq:cumulant}). It can be shown that both
definitions are equivalent for all cumulants of order greater than one
because the classical cumulants are translation invariant. The
definition in (\ref{eq:cumulant}) is thus essentially the same as~(\ref{eq:kappa-in-rho}) but with a different defining poset (cf. \cite
{rotacumulants,speed1983cumulants}).

Using a basic result in the theory of lattices, Lemma \ref
{lem:splitzero} shows that certain features of classical cumulants are
also shared by tree cumulants (cf. Section 2.1 of
\cite{mccullagh1987tms}).\vspace*{-2pt}
\begin{lem}[(Corollary in \cite{rota1964fct}, Section 5)]\label{lem:stanley}
Let $L$ be a finite lattice and let $\pi_{0}\neq\hat{1}$ in $L$.
Then, for any $\nu$ in $L$
\[
\sum_{\pi\wedge\pi_{0}=\nu}\mm(\pi,\hat{1})=0.\vspace*{-2pt}
\]
\end{lem}

\begin{lem}\label{lem:splitzero}
Let $T$ be a tree with $n$ leaves. Whenever there exists an edge split
$C_{1}|C_{2}\in\Pi_{T}$ of the set of leaves $[n]$ such that
$X_{C_{1}}\indep X_{C_{2}}$, then $\kappa_{1\cdots n}=0$.\vspace*{-2pt}
\end{lem}
\begin{pf}
Let $\pi_{0}$ be the split $C_{1}|C_{2}$ such that $X_{C_{1}}\indep
X_{C_{2}}$. It follows that $\mu_{1\cdots n}$ is equal to $\mu
_{C_{1}}\mu_{C_{2}}$. More generally, for any $I\in[n]_{\geq2}$,
\[
\mu_{I}=\mu_{C_{1}\cap I}\mu_{C_{2}\cap I}.\vspace*{-2pt}
\]
Consequently, for any partition $\pi\in\Pi_{T}$
\begin{equation}\label{eq:factmomcum}
\prod_{B\in\pi}\mu_{B}= \prod_{B\in\pi\wedge\pi_{0}} \mu_{B}.\vspace*{-2pt}
\end{equation}
Using (\ref{eq:kappa-in-rho}) and (\ref{eq:factmomcum}), we obtain
\[
\kappa_{1\cdots n}=\sum_{\pi\in\Pi_{T}}\mm(\pi,\hat{1})\prod
_{B\in\pi} \mu_{B}=\sum_{\pi\in\Pi_{T}}\mm(\pi,\hat{1})\prod
_{B\in\pi\wedge\pi_{0}} \mu_{B}.\vspace*{-2pt}\vadjust{\goodbreak}
\]
Since $\pi\wedge\pi_{0}\leq\pi_{0}$, by grouping all partitions
$\pi\in\Pi_{T}$ giving the same partition, after taking the meet
with $\pi_{0}$, we can rewrite the sum as
\[
\kappa_{1\cdots n}=\sum_{\pi\in\Pi_{T}}\mm(\pi,\hat{1})\prod
_{B\in\pi\wedge\pi_{0}} \mu_{B}=\sum_{\nu\leq\pi_{0}}
\biggl(\sum_{\pi\wedge\pi_{0}=\nu}\mm(\pi,\hat{1}) \biggr)\prod
_{B\in\pi\wedge\pi_{0}} \mu_{B}.
\]
However, this is zero since by Lemma \ref{lem:stanley} each of $\sum
_{\pi\wedge\pi_{0}=\nu}\mm(\pi,\hat{1})$ is zero.
\end{pf}

\section{The induced parametrization}\label{sec:parametrization}

We now define a new parameter space, $\Omega_{T}$, with $|V|+|E|$
parameters denoted by $\eta_{u,v}$ for all $(u,v)\in E$ and $\bar{\mu
}_{v}$ for all $v\in V$. The map between the two parameter spaces is
given by
\begin{eqnarray}\label{eq:uij}
\eta_{u,v}&=&\theta^{(v)}_{1|1}-\theta^{(v)}_{1|0} \qquad\mbox{for
all $(u,v)\in E$}  \quad \mbox{and}\nonumber
\\[-8pt]
\\[-8pt]
 \bar{\mu}_v&=&1-2\lambda_v\qquad\mbox{for each } v\in V,
\nonumber
\end{eqnarray}
where $\lambda_{v}$ is a polynomial in the original parameters in
$\Theta_{T}$. The details are given in Appendix~\ref
{sec:repar-models}, where the inverse map is given by (\ref
{eq:theta-inverse}). It follows that the change of parameters between
$\Theta_{T}$ and $\Omega_{T}$ is a polynomial isomorphism.

It can be checked that if $\operatorname{{Var}}(Y_{u})>0,$ then $\eta_{u,v}=\E
(U_{u}U_{v})/\operatorname{{Var}}( Y_u)$ is the regression coefficient of $Y_{v}$
on $Y_{u}$. Therefore, $\eta_{u,v}$, defined above, coincides with the
definition of~$\eta_{u,v}$ in (\ref{eq:eta-dajesz}). If $\operatorname{{Var}}(Y_{u})=0,$ then the formula in (\ref{eq:eta-dajesz}) is not well
defined; however, (\ref{eq:uij}) always is.


Proposition \ref{prop:monomial} below motivates the whole section and
demonstrates why our new coordinate system is particularly useful.
Henceforth let $\cM_T^\kappa=(f_{\mu\kappa}\circ f_{p\mu})(\cM
_T)\subseteq\cK_T$.
\begin{prop}\label{prop:monomial}
Let $T=(V,E)$ be a rooted tree with $n$ leaves such that each inner
node has degree at most three. Then $\cM_T^\kappa$ is given as the
image of $\psi_{T}\dvtx\Omega_{T}\rightarrow\cK_{T}$. Here~$\psi_{T}$
is defined by $\lambda_{i}=\frac{1}{2}(1-\bar{\mu}_{i})$ for
$i=1,\ldots,n$ and
\begin{equation}\label{eq:kappa_def_general}
\kappa_{I}=\frac{1}{4} \bigl(1-\bar{\mu}_{r(I)}^{2} \bigr) \prod
_{v\in V(I)\setminus I} \bar{\mu}_{v}^{\deg(v)-2}\prod_{(u,v)\in
E(I)} \eta_{u,v}\qquad\mbox{for each } I\in[n]_{\geq2},
\end{equation}
where the degree is taken in $T(I)=(V(I),E(I))$ and $r(I)$ denotes the
root of $T(I)$ (cf. Section~\ref{sec:intro_trees}).
\end{prop}

The proof is given in Appendix \ref{app:prop}.

By Lemma \ref{lem:triv-red} we can obtain the parametrization of $\cM
_{T}$ for any non-trivalent tree $T=(V,E)$ using a parametrization for
its trivalent expansion $T^{*}=(V^{*},E^{*})$. Let $E'$ be the subset
of inner nodes of $E^{*}$ given in Definition \ref{def:triv-exp}, so
that ${T}^{*}/E'=T$. Let $\{V^{*}\}$ denote the equivalence classes of
subsets of $V^{*}$ such that $v\sim v'$ if and only if $v$
becomes
identified with $v'$ in $T$ in the process of contracting $E'$ in
$T^{*}$. There exists a natural identification of $V$ with $\{V^{*}\}$.
Let $\{v\}$ denote the equivalence class of $v\in V^{*}$ or the
corresponding node in $T$. In particular, since $E'$ is a set of inner
edges, the class $\{i\}$ of every leaf $i\in[n]$ can be naturally
identified with $i$ and hence $\{V^{*}\setminus[n]\}=\{V^{*}\}
\setminus[n]$.
\begin{lem} Let $T$ be any tree and ${T}^{*}$ be its trivalent
expansion. If $\kappa_{I}^{*}$ for $I\in[n]_{\geq2}$ are tree
cumulants of ${T}^{*}$, then $\cM_{T}^{\kappa}$ is given in $\cK
_{T^{*}}$ as the image of a map that is the identity on the coordinates
corresponding to $\bar{\mu}_{i}$ for $i=1,\ldots,n$ and, for each
$I\in[n]_{\geq2}$,
\begin{equation}\label{eq:forcontracted}
\kappa_{I}^{*}=\frac{1}{4} \bigl(1-\bar{\mu}_{r(I)}^{2} \bigr)
\prod_{v\in V(I)\setminus I} \bar{\mu}_{v}^{\deg(v)-2}\prod
_{(u,v)\in E(I)} \eta_{u,v},
\end{equation}
where $T(I)=(V(I),E(I))$ is the subtree of $T$ spanned on $I$.
\end{lem}
\begin{pf}
By Lemma \ref{lem:triv-red} and equation (\ref{eq:star-constraints}),
$\cM_{T}\subseteq\cM_{{T}^{*}}$ is the image $f_{{T}^{*}}(\Theta
_{T})$, where~$\Theta_{T}$ is the subset of $\Theta_{{T}^{*}}$ given
by setting ${\theta^{*}}^{(v)}_{\a_{v}|\a_{u}}=\delta_{\a_{u}\a
_{v}}$ for every edge\vspace*{2pt} $(u,v)\in E'$ and ${\theta^{*}}^{(v)}_{\a
_{v}|\a_{u}}={\theta}^{(v)}_{\a_{v}|\a_{u}}$ otherwise. In the new
parameters, $\Omega_{T}$ is isomorphic to the subset of~$\Omega
_{{T}^{*}}$ given by
\begin{eqnarray}\label{eq:constr-omega}
\eta_{u,v}^{*}&=&\eta_{u,v}\qquad\mbox{for all } (u,v)\notin E',\nonumber\\
\eta_{u,v}^{*}&=&1 \qquad\mbox{for all } (u,v)\in E'  \quad \mbox{and}\\
\bar{\mu}_{v}^{*}&=&\bar{\mu}_{\{v\}}\qquad\mbox{for all } v\in V^{*}.\nonumber
\end{eqnarray}
Denote the root of $T^{*}$ by $r^{*}$. We show (\ref
{eq:forcontracted}) for $I=[n]$. The general case can be proved with an
obvious change in notation. By Proposition \ref{prop:monomial}, the
model $\cM_{{T}^{*}}$ is parametrized~by
\begin{equation}\label{eq:kappa-star}
\kappa_{1\cdots n}^{*}=\frac{1}{4} (1-{\bar{\mu
}}_{r^{*}}^{*2} ) \prod_{v\in{V}^{*}\setminus[n]} \bar{\mu
}_{v}^{*\deg(v)-2}\prod_{(u,v)\in{E}^{*}} \eta_{u,v}^{*}.
\end{equation}
Since ${E}^{*}=E\cup E'$ by applying (\ref{eq:constr-omega}), $\prod
_{(u,v)\in{E}^{*}} \eta_{u,v}^{*}$ becomes $\prod_{(u,v)\in{E}}
\eta_{u,v}$, where we have identified $E$ with $E^{*}\setminus E'$.
For every $w\in V^{*}$, whenever $\deg\{w\}\geq3$, we have that $\deg
\{w\}=|\{w\}|+2$. Therefore, if $\deg\{w\}\geq3,$ then the degree of
each $v\in\{w\}$ in $T^{*}$ equals $3$. Hence
\[
\sum_{v\in\{w\}}(\deg v-2)=\sum_{v\in\{w\}} 1=|\{w\}|=\deg\{w\}-2.
\]
It follows that, after applying (\ref{eq:constr-omega}), $\prod_{v\in
\{w\}}\bar{\mu}^{*\deg v-2}_{v}$ becomes $\bar{\mu}_{\{w\}}^{\deg\{
w\} -2}$. The last statement is also true if $\deg\{w\}=2$. For, in
this case, $\deg w=2$ in $T^{*}$ and $w$ is the only element in $\{w\}
$. Moreover, $E'$ is necessarily contained in the set of inner edges of
$T^{*}$. It follows that $\prod_{v\in{V}^{*}\setminus[n]} \bar{\mu
}_{v}^{*\deg(v)-2}$ in (\ref{eq:kappa-star}) becomes
\[
\prod_{\{w\}\in\{{V^{*}}\}\setminus[n]} \bar{\mu}_{\{w\}}^{\deg(\{
w\})-2}=\prod_{v\in{V}\setminus[n]} \bar{\mu}_{v}^{\deg(v)-2}.\vadjust{\goodbreak}
\]
In addition, $\{r^{*}\}$ becomes the root of $T$ denoted by $r$.
Therefore, (\ref{eq:kappa-star}) becomes
\[
\kappa_{1\cdots n}^{*}=\frac{1}{4} (1-\bar{\mu}_{r}^{2}
) \prod_{v\in V\setminus[n]} \bar{\mu}_{v}^{\deg(v)-2}\prod
_{(u,v)\in E} \eta_{u,v},\vspace*{-2pt}
\]
which is exactly (\ref{eq:forcontracted}) for $I=[n]$.\vspace*{-2pt}
\end{pf}
%

\begin{rem}\label{rem:phyl}
For every $v\in V$ the variance $\operatorname{{Var}}(Y_{v})$ is zero if and only
if $\bar{\mu}_{v}^{2}=1$. Hence, in the case when $\bar{\mu
}_{v}^{2}<1,$ the variable $Y_{v}$ is non-degenerate. In phylogenetics
it is usually assumed that $\bar{\mu}_{{r}}^{2}<1$ for the root $r$
of $T$ and $\eta_{u,v}\neq0$ for all $(u,v)\in E$ (cf.
Conditions~(M1) and (M2) in Section 8.2, \cite{semple2003pol}). It is shown in
Section 8.2 in \cite{semple2003pol} that~(M1) and~(M2) imply the
weaker condition $\bar{\mu}_{{v}}^{2}<1$ for all $v\in V$. Over the
subset of $\Omega_{T}$ on which this weaker condition holds, we can
apply another smooth transformation on both the parameter and model
space. This leads to a further simplification of the parametrization in
(\ref{eq:kappa_def_general}) presented in Appendix \ref{app:nondeg}.\vspace*{-3pt}
\end{rem}


\section{Singularities and the geometry of unidentified
subspaces}\vspace*{-3pt}\label{sec:fibers}

The identifiability of general Markov models can be addressed here
geometrically. For any $q\in\cM_{T}$ the preimage $\widehat{\Theta
}_{T}:=f_{T}^{-1}(q)$, that is, the set of parameter values that is
consistent with the known probability model $q$, is called the
$q$\textit{-fiber}. In this section, we analyze the geometry of these
fibers, determining when they are finite and thus when the model is
locally identifiable. We will also be interested in when the fibers are
smooth subsets of $\Theta_{T}$ and when they are singular. We use
methods similar to the ones presented in a different context by Moulton
and Steel in \cite{moulton2004ppo}, Section 6. The results in this
section generalize similar results for the naive Bayes models (cf.
\cite{geiger2001sef}, Theorem 7).

First we analyze the geometric description of $\Omega_{T}$. This gives
a set of implicit inequalities constraining each $q$-fiber. Simple
linear constraints defining $\Theta_{T}$ become only slightly more
complicated when expressed in the new parameters. The choice of
parameter values is not free anymore in the sense that the constraining
equations for each of the parameters involve the values of other
parameters. By (\ref{eq:theta-inverse}), $\Omega_{T}$ is given by
$\bar{\mu}_{r}\in[-1,1]$ and for each $(u,v)\in E$
\begin{eqnarray}\label{eq:constraints}
-(1+\bar{\mu}_{v}) &\leq&(1-\bar{\mu}_{u})\eta_{u,v} \leq(1-\bar
{\mu}_{v}),\nonumber
\\[-9pt]
\\[-9pt]
-(1-\bar{\mu}_{v}) &\leq&(1+\bar{\mu}_{u})\eta_{u,v} \leq(1+\bar
{\mu}_{v}).
\nonumber\vspace*{-2pt}
\end{eqnarray}

For $\hat{p}\in\cM_T$ let $\widehat{\Sigma}=[\hat{\mu}_{ij}]\in
\R^{n\times n}$ be the covariance matrix of the observed variables
labelled by the leaves of $T$ computed with respect to $\hat{p}$. We
show that the geometry of the $\hat{p}$-fiber, denoted by $\widehat
{\Theta}_{T}$, is determined by zeros in $\widehat{\Sigma}$. Let
$\hat{\lambda}_{i}$ be the expected value of $X_{i}$. Then, for every
point in the $\hat{p}$-fiber, we have $\bar{\mu}_{i}=\hat{\mu
}_{i}=1-2\hat{\lambda}_{i}$ for all $i=1,\ldots,n$. Without loss we
always assume that $\hat{\lambda}_{i}(1-\hat{\lambda}_{i})\neq0$
(or, equivalently, that $\hat{\mu}_{i}^{2}\neq1$) for all
$i=1,\ldots,n$.

It is easier to analyze the geometry of $\hat{p}$-fibers in $\Omega
_{T}$. Therefore transform $\widehat{\Theta}$ to $\Omega_{T}$ using
the mapping\vadjust{\goodbreak} $f_{\theta\omega}$. The image of this map, denoted by
$\widehat{\Omega}_{T}$, is isomorphic to $\widehat{\Theta}_{T}$.
Let $\hat{\kappa}_{ij}$ denote the corresponding second-order tree
cumulants in the point $f_{p\kappa}(\hat{p})$. Since $\kappa
_{ij}=\mu_{ij}$ for all $i,j\in[n]$, from (\ref
{eq:kappa_def_general}) for any $\omega_{0}=((\bar{\mu
}_{v}^{0}),(\eta_{u,v}^{0}))\in\widehat{\Omega}_{T}$ we have that
\begin{equation}\label{eq:muij}
\hat{\mu}_{ij}=\mu_{ij}(\omega_{0})=\frac{1}{4}\bigl (1-\bigl({\bar
{\mu}^{0}_{r(ij)}}\bigr)^2 \bigr) \prod_{(u,v)\in E(ij)}
\eta_{u,v}^{0}.\vspace*{-2pt}
\end{equation}
We say that that an edge, $e\in E,$ is \textit{isolated relative to}
$\hat{p}$ if $\hat{\mu}_{ij}=0$ for all\vspace*{1pt} $i,j\in[n]$ such that $e\in
E(ij)$. We denote the set of all edges of $T$ that are isolated
relative to $\hat{p}$ by $\widehat{E}\subseteq E$. We define the
$\hat{p}$-forest $\widehat{T}$ as the forest obtained from $T$ by
removing edges in $\widehat{E}$ so that $\widehat{T}=T\setminus
\widehat{E}$. Hence, the set of vertices of $\widehat{T}$ is equal to
the set of vertices of $T$ and the set of edges is equal to $E\setminus
\widehat{E}$.

\begin{figure}
\centerline{
\begin{minipage}{134pt}
\begin{figure}[H]
\includegraphics{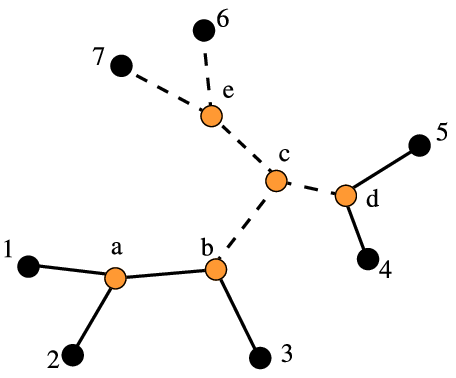}
\end{figure}
\end{minipage}\qquad
\begin{minipage}{105pt}
\begin{figure}[H]
$\widehat{\Sigma}=
\left[
%
\begin{array}{ccccccccccccc}* && * && * && 0 && 0 && 0 && 0 \\ && * && * && 0
&& 0 && 0 && 0 \\ && && * && 0 && 0 && 0 && 0 \\ && && && * && * && 0
&& 0 \\ && && && && * && 0 && 0 \\ && && && && && * && 0 \\
&& && && && && && *
\end{array}
 \right]$
\end{figure}
\end{minipage}
}\vspace*{-3pt}
\caption{An example of a tree and a sample covariance matrix. The
dashed lines depict the edges isolated with respect to $\hat
{p}$.}\label{fig:isolated}
\vspace*{-4pt}
\end{figure}

We illustrate this construction in the example below. Let $T$ be the
tree given in Figure~\ref{fig:isolated} and assume that the covariance
matrix contains zeros given in the provided $7\times7$ matrix, where
the asterisks mean any non-zero values such that the matrix is positive
semidefinite.
It can be checked that $\widehat{E}=\{(b,c),(c,d),(c,e),(e,6), (e,7)\}
$ and these edges are depicted as dashed lines. The forest, $\widehat
{T}$, is obtained by removing the edges in $\widehat{E}$.\vspace*{2pt}

%

We now define relations on $\widehat{E}$ and $E\setminus\widehat
{E}$. For two edges, $e,e'$, with either $\{e,e'\}\subset\widehat{E}$
or $\{e,e'\}\subset E\setminus\widehat{E}$, write $e\sim e'$ if
either $e=e'$ or $e$ and $e'$ are adjacent and all the edges that are
incident with both $e$ and $e'$ are isolated relative to $\hat{p}$. We
now construct the transitive closure of $\sim$ restricted to pairs of
edges in $\widehat{E}$ to form an equivalence relation on $\widehat
{E}$. Consider a graph with nodes representing elements of $\widehat
{E}$ and put an edge between~$e,e'$ whenever $e\sim e'$. Then the
equivalence classes correspond to connected components of this graph.
In the same way, we take the transitive closure of $\sim$ restricted
to the pairs of edges in $E\setminus\widehat{E}$ to form an
equivalence relation in $E\setminus\widehat{E}$. We will let
$[\widehat{E}]$ and $[E\setminus\widehat{E}]$ denote the set of
equivalence classes of $\widehat{E}$ and $E\setminus\widehat{E}$,
respectively. For the tree from the example above, $[\widehat{E}]$ is
one element given by a subtree of $T$ spanned on $\{b,d,6,7\}$
and
\[
[E\setminus\widehat{E}]= \{\{(1,a)\}, \{(2,a)\}, \{(a,b), (b,3)\}
, \{(d,4), (d,5)\}  \}.\vspace*{-2pt}
\]

By construction, all the inner nodes of $T$ have either degree zero in
$\widehat{T}$ or the degree is strictly greater than one. The
following lemma shows that whenever the degree of an inner node in
$\widehat{T}$ is not zero, the node represents a non-degenerate random
variable.\vadjust{\goodbreak}

\begin{lem}\label{lem:non-deg}
Let $\hat{p}\in\cM_{T}$. If $v\in V$ is an inner node of $T$ such
that $\deg(v)\geq2$ in the $\hat{p}$-forest $\widehat{T}$, then the
variable $H_{v}$ cannot be degenerate.
\end{lem}
\begin{pf} By construction, if $\deg(v)\geq2$ in $\widehat{T}$,
then there exists $i,j\in[n]$ such that $\hat{\mu}_{ij}\neq0$ and
$v$ lies on the path between $i$ and $j$. Suppose that $H_{v}$ is
degenerate. Then the global Markov properties in (\ref{eq:GMP}) imply
that $X_{i}\indep X_{j}$. But then $\hat{\mu}_{ij}=0$ and we obtain
the contradiction.
\end{pf}

We now list some basic statements, partly based on Lemma 6.4 in \cite
{moulton2004ppo}, which follow directly definitions above.
\begin{rem}\label{rem:steel_moulton}
Let $T=(V,E)$ be a tree with $n$ leaves, let $\cM_{T}$ be the
corresponding general Markov model and suppose that $\hat{p}\in\cM_T$.
\begin{longlist}[(iii)]
\item[(i)] The edges in any equivalence class of $[\widehat{E}]$ form
a connected subgraph of $T$. If $T$ is trivalent, then this subgraph is
either a single edge or a trivalent tree.
\item[(ii)] If each inner node of $T$ has degree at least two in
$\widehat{T}$, then all the equivalence classes in $[\widehat{E}]$
are just single edges. If each inner node has degree at least three in
$\widehat{T}$, then\vspace*{2pt} all equivalence classes in $[E\setminus\widehat
{E}]$ are single edges.
\item[(iii)] The edges in any equivalence class in $[E\setminus
\widehat{E}]$ can be ordered so that they form a~path in~$T$.
\item[(iv)] Every connected component of $\widehat{T}$ is either a
single node or a tree with its set of leaves contained in $[n]$.
\end{longlist}
\end{rem}

\begin{lem}\label{lem:steel_moulton}
Let $E(uv)\subset E$ be any path as in Remark \ref{rem:steel_moulton}\textup{(iii)},
which is an element of $[E\setminus\widehat{E}]$. Then the
quantities $\mu_{uv}^{2}$ and $\eta_{u,v}^{2}$ are constant on
$\widehat{\Omega}_{T}$ and non-zero. It is possible to determine
their values from $\hat{p}$.
\end{lem}
\begin{pf}
First note that the degree of each inner node on the path between $u$
and $v$ in $\widehat{T}$ must be exactly two. Moreover, the degree of
both $u$ and $v$ in $\widehat{T}$ must be at least three unless $u$ or
$v$ is a leaf. Consider the case when both $u$ and $v$ are inner nodes
of $T$. In this case, these nodes have degrees at least three in
$\widehat{T}$ and we can find four leaves $i,j,k,l$ such that $u$
separates $i$ from $j$ in $\widehat{T}$, $v$ separates $k$ and $l$ and
$\{u,v\}$ separates $\{i,j\}$ from $\{k,l\}$ as in the graph below.
\begin{center}\vspace*{6pt}

\includegraphics{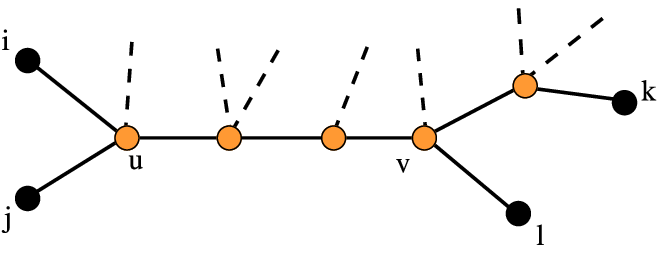}
\vspace*{6pt}
\end{center}
Furthermore, by construction, $\hat{\mu}_{ij},\hat{\mu}_{kl},\hat
{\mu}_{ik},\hat{\mu}_{jl}$ are all non-zero. Consider the marginal
models for $T(ijk)$ and $T(ikl)$. By Corollary \ref{cor:marg-trip},
these are equivalent to models associated with tripod trees as in
Figure \ref{fig:tripod}. Hence, from (\ref{eq:params-tripod}) we have that
\begin{equation}\label{eq:means}
\bar{\mu}_{u}^{2}=\frac{\hat{\mu}_{ijk}^{2}}{\hat{\mu
}_{ijk}^{2}+4\hat{\mu}_{ij}\hat{\mu}_{ik}\hat{\mu}_{jk}},\qquad
\bar{\mu}_{v}^{2}=\frac{\hat{\mu}_{ikl}^{2}}{\hat{\mu
}_{ikl}^{2}+4\hat{\mu}_{ik}\hat{\mu}_{il}\hat{\mu}_{kl}}.
\end{equation}
These equations are well defined since $\hat{\mu}_{ij}\hat{\mu
}_{ik}\hat{\mu}_{jk}>0$ and $ \hat{\mu}_{ik}\hat{\mu}_{il}\hat
{\mu}_{kl}>0$.
Consider the quantity $\frac{\hat{\mu}_{ik}\hat{\mu}_{jl}}{\hat
{\mu}_{ij}\hat{\mu}_{kl}}$ and substitute (\ref{eq:muij}) for each
of the terms. A simple rearrangement now gives that
\[
\frac{\hat{\mu}_{ik}\hat{\mu}_{jl}}{\hat{\mu}_{ij}\hat{\mu
}_{kl}}=\frac{1-\bar{\mu}_{u}^{2}}{1-\bar{\mu}_{v}^{2}}\eta
_{u,v}^{2}(\omega),
\]
where $\eta_{u,v}(\omega)=\frac{1-\bar{\mu}_{r(uv)}^{2}}{1-\bar
{\mu}_{u}^{2}}\prod_{(w,w')\in E(uv)} \eta_{w,w'}$. Therefore,
substituting for $\bar{\mu}_{u}^{2}$, $\bar{\mu}_{v}^{2}$ using~(\ref{eq:means}) implies that $\eta_{u,v}^{2}$ is constant on
$\widehat{\Omega}_{T}$ and non-zero. Its value can be determined as a~function of $\hat{p}$. Also the value of $\mu_{uv}^{2}$ is constant
since $\mu_{uv}^{2}=\frac{1}{16}(1-\bar{\mu}_{u}^{2})^{2}\eta^{2}_{u,v}$.

If either $u$ or $v$ is a leaf of $T$, then the argument is very
similar. Thus, if $u$ is a leaf, then consider any two leaves $i,j$ of
$T$ such that $v$ separates $u,i,j$ in $\widehat{T}$. In particular,
as in (\ref{eq:means}),
\[
 \bar{\mu}_{v}^{2}=\frac{\hat{\mu}_{uij}^{2}}{\hat{\mu
}_{uij}^{2}+4\hat{\mu}_{ui}\hat{\mu}_{uj}\hat{\mu}_{ij}}.
\]
Moreover, $\eta_{u,v}(\omega)$ must be determined, since from (\ref{eq:muij})
\[
\frac{\hat{\mu}_{ui}\hat{\mu}_{uj}}{\hat{\mu}_{ij}}=\frac
{1}{4}(1-\bar{\mu}_{v}^{2})\eta_{u,v}^{2}(\omega),
\]
from which it follows that $\eta^{2}_{u,v}$ has to be constant on the
$\hat{p}$-fiber.
\end{pf}

The following theorem shows that the geometry of the $\hat{p}$-fiber
$\widehat{\Omega}_T$ is determined by the zeros of the covariance
matrix $\widehat{\Sigma}$.
\begin{theo}[(The geometry of the $\bolds{\hat{p}}$-fiber -- the smooth case)]
\label{lem:fibers}Let $\hat{p}\in\cM_{T}$. If each of the inner
nodes of $T$ has degree at least three in the $\hat{p}$-forest
$\widehat{T}$, then the $\hat{p}$-fiber is a finite set of points of
cardinality $2^{|V|-n}$. If each of the inner nodes of $T$ has degree
at least two in $\widehat{T}$, then the $\hat{p}$-fiber is
diffeomorphic to a disjoint union of polyhedra. In particular, it is a
manifold with corners. Its dimension is $2l_{2}$, where $l_{2}$ is the
number of degree-$2$ nodes in $\widehat{T}$.
\end{theo}

The proof is given in Appendix \ref{app:proofs}.

If $T$ is trivalent, then the $\hat{p}$-fiber is finite if and only if
for all $i,j\in[n]$ $\mu_{ij}\neq0$. The proof of Theorem \ref
{lem:fibers} provides explicit formulae for the parameters in this case
when the $\hat{p}$-fiber is a finite number of points.
\begin{cor}\label{cor:formulas}
Let $T$ be a tree such that each inner node has degree at least three
and let $\hat{p}\in\cM_{T}$. Consider the\vadjust{\goodbreak} $\hat{p}$-forest
$\widehat{T}$. If every inner node of $T$ has degree at least three in~%
$\widehat{T}$, then, by Remark \ref{rem:steel_moulton}\textup{(ii)}, both
$[\widehat{E}]$ and $[E\setminus\widehat{E}]$ consist of singletons.
In this case, every point in the $\hat{p}$-fiber satisfies
\begin{eqnarray}
\bar{\mu}_{i}&=&\hat{\mu}_{i}  \qquad  \mbox{for all }
i=1,\ldots,n,\nonumber
\\[-10pt]
\\[-10pt]
\eta_{u,v}&=&0  \qquad  \mbox{for all } (u,v)\in[\widehat{E}].
\nonumber\vspace*{-2pt}
\end{eqnarray}
Moreover, for any inner node $v$ of $T$, if $i,j,k\in[n]$ are any
three leaves separated by $v$ in~$T$ such that $\hat{\mu}_{ij}\hat
{\mu}_{ik}\hat{\mu}_{jk}\neq0,$ then
\[
\bar{\mu}_{v}^{2}=\frac{\hat{\mu}_{ijk}^{2}}{\hat{\mu
}_{ijk}^{2}+4\hat{\mu}_{ij}\hat{\mu}_{ik}\hat{\mu}_{jk}}\vspace*{-2pt}
\]
for any terminal edge $(v,i)\in E\setminus\widehat{E}$, where $v$ is
an inner node and $i\in[n]$ is a leaf of $T$. Let $j,k$ be any two
leaves such that $v$ separates $i,j,k$ and $\hat{\mu}_{jk}\neq0$. Then
\[
\eta_{v,i}^{2}=\frac{\hat{\mu}_{ijk}^{2}+4\hat{\mu}_{ij}\hat{\mu
}_{ik}\hat{\mu}_{jk}}{\hat{\mu}^{2}_{jk}}.\vspace*{-2pt}
\]
Moreover, for any inner edge $(u,v)\in E\setminus\widehat{E}$ let
$i,j,k,l\in[n]$ be any four leaves of $T$ such that $u$ separates $i$
and $j$ in $\widehat{T}$, $v$ separates $j$ and $k$ in $\widehat{T}$
and $(u,v)$ separates $\{i,j\}$ from $\{k,l\}$ in $\widehat{T}$. Then
\[
\eta_{u,v}^{2}=\frac{\hat{\mu}_{il}^{2}}{\hat{\mu}_{ij}^{2}}\frac
{\hat{\mu}_{ijk}^{2}+4\hat{\mu}_{ij}\hat{\mu}_{ik}\hat{\mu
}_{jk}}{\hat{\mu}_{ikl}^{2}+4\hat{\mu}_{ik}\hat{\mu}_{il}\hat
{\mu}_{kl}}.\vspace*{-2pt}
\]
\end{cor}

\begin{rem}
The choice of signs of the $\bar{\mu}_{v}$ and $\eta_{u,v}$ in
Corollary \ref{cor:formulas} is not completely free and has to be
consistent with signs of tree cumulants via (\ref
{eq:kappa_def_general}) (see Appendix \ref{app:signs}).\vspace*{-2pt}
\end{rem}

The singular case when there is at least one degree-zero inner node is
more complicated. We begin with an example.\vspace*{-2pt}

\begin{exmp}\label{ex:sing}
Let $T=(V,E)$ be the tripod tree rooted in the inner node as in Figure
\ref{fig:tripod} and let $\hat{p}\in\cM_{T}$. The degree of $h$ in
the $\hat{p}$-forest $\widehat{T}$ is less than two if and only if
$\hat{\mu}_{ij}=0$ for all $i\neq j=1,2,3$. In this situation,
$\widehat{E}=E$ and the $\hat{p}$-fiber $\widehat{\Omega}_T$ is
given as a subset of $\Omega_{T}$ by equations for the sample means
$\bar{\mu}_{i}=\hat{\mu}_{i}$ for $i=1,2,3$ together with the three
additional equations
\[
(1-\bar{\mu}_h^2) \eta_{h,1} \eta_{h,2}=0, \qquad  (1-\bar{\mu}_h^2)
\eta_{h,1} \eta_{h,3}=0,  \qquad  (1-\bar{\mu}_h^2) \eta_{h,2} \eta_{h,3}=0.\vspace*{-2pt}
\]
Geometrically, in the subspace given by $\bar{\mu}_{i}=\hat{\mu
}_{i}$ for $i=1,2,3$, this is a union of two three-dimensional
hyperplanes $\{\bar{\mu}_h=\pm1\}$ and three planes given by $\{\eta
_{h,1}=\eta_{h,2}=0\}$, $\{\eta_{h,1}=\eta_{h,3}=0\}$ and $\{\eta
_{h,2}=\eta_{h,3}=0\}$ subject to the additional inequality
constraints defining $\Omega_{T}$ and given by (\ref
{eq:constraints}). In particular, it is not a regular set since it has
self-intersection points given by $1-\bar{\mu}_{h}^{2}=\eta
_{h,1}=\eta_{h,2}=\eta_{h,3}=0$.\vspace*{-2pt}
\end{exmp}


This geometry is mirrored in the general case. We first need two
definitions. We say that a node $v\in V$ is\vadjust{\goodbreak} \textit{non-degenerate
(with respect to $\hat{p}$)} if either $v$ is a leaf of $T$ or $\deg
v\geq2$ in $\widehat{T}$. Otherwise, we say that the node is \textit
{degenerate with respect to} $\hat{p}$. The set of all nodes that are
degenerate with respect to $\hat{p}$ is denoted by $\widehat{V}$. By
Lemma \ref{lem:non-deg}, for all $v\in V\setminus\widehat{V}$, $\operatorname{{Var}}(Y_{v})\neq0$, where the variance is computed with respect to
$\hat{p}$. Hence $v$ is non-degenerate if and only if $Y_{v}$ is a
non-degenerate random variable.

We define the \textit{deepest singularity} of $\widehat{\Omega}_{T}$ as
\begin{equation}\label{eq:deepest}
\widehat{\Omega}_{\mathrm{deep}}:=\{\omega\in\widehat{\Omega}_T\dvt
\eta_{u,v}=0,  \bar{\mu}_{v}^{2}=1 \mbox{ for all } (u,v)\in
\widehat{E},  v\in\widehat{V}\}.\vspace*{-3pt}
\end{equation}
%

\begin{theo}[(The geometry of the $\bolds{\hat{p}}$-fiber -- the singular
case)]\label{prop:singular}
If $\widehat{V}$ is non-empty, then the $\hat{p}$-fiber is a singular
variety given as a union of intersecting smooth manifolds in $\R
^{|V|+|E|}$ restricted to $\Omega_{T}$. Their common intersection
locus restricted to $\Omega_{T}$ is given by $\widehat{\Omega}_{\mathrm{deep}}$,
which lies on the boundary of $\Omega_{T}$.\vspace*{-3pt}
\end{theo}

The proof is given in Appendix \ref{app:proofs}.

\begin{figure}

\includegraphics{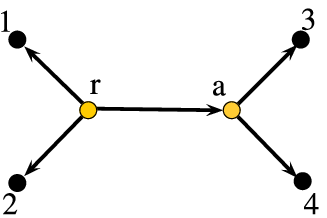}

\caption{The quartet tree.}\label{fig:quartet}
\end{figure}

\vspace*{-3pt}
\section{Example: The quartet tree model}\label{sec:quartet}
\vspace*{-3pt}

In this section, we study the first non-trivial example: the quartet
tree model given by the tree in Figure \ref{fig:quartet}. The model is
parametrized as in (\ref{eq:p_albar2}) by the root distribution and
conditional probabilities attached to each of the edges. We set the
values of the parameters to $\theta^{(r)}_{1}=0.8$, $\theta
^{(1)}_{1|0}=0.8$, $\theta^{(1)}_{1|1}=0.3$, $\theta
^{(2)}_{1|0}=0.7$, $\theta^{(2)}_{1|1}=0.3$, $\theta
^{(a)}_{1|0}=0.8$, $\theta^{(a)}_{1|1}=0.3$, $\theta
^{(3)}_{1|0}=0.7$, $\theta^{(3)}_{1|1}=0.3$, $\theta
^{(4)}_{1|0}=0.7$, $\theta^{(4)}_{1|1}=0.3$.
Using (\ref{eq:p_albar2}) we can then calculate the corresponding
probabilities over the observed nodes that are given in the third
column in the table below. The change of coordinates $f_{p\lambda}$
presented in Appendix \ref{sec:change} and $f_{\mu\kappa}$ in
Section \ref{sec:kappas} gives the corresponding non-central moments
and tree cumulants that are shown in Table \ref{tab:quartet1}.
Formula (\ref{eq:uij}) enables us to calculate the values for the new
parameters as: $\eta_{r,1}=0.5$, $\eta_{r,2}=0.4$, $\eta_{r,a}=0.5$,
$\eta_{a,3}=0.4$, $\eta_{a,4}=0.4$ and $\bar{\mu}_{1}=-0.4$, $\bar
{\mu}_{2}=-0.24$, $\bar{\mu}_{3}=-0.16$, $\bar{\mu}_{4}=-0.16$,
$\bar{\mu}_{r}=-0.6$, $\bar{\mu}_{a}=-0.4$. It is easy to verify
that (\ref{eq:kappa_def_general}) holds in this example. For instance,
\[
\kappa_{1234}=\tfrac{1}{4}(1-\bar{\mu}_{r}^{2})\bar{\mu}_{r}\bar
{\mu}_{a}\eta_{r,1}\eta_{r,2}\eta_{r,a}\eta_{a,3}\eta_{a,4}=0.0006,
\]
which equates with the value in the table. In general, higher-order
tree cumulants tend to be very small.

\begin{table}
\tablewidth=168pt
\caption{Moments and tree cumulants for a probability assignment in
$\cM_{T}$}\label{tab:quartet1}
\begin{tabular}{@{}ld{4.0}lld{2.4}@{}}
  \hline
$\a$& \multicolumn{1}{l}{$I$} & $p_{\a}$ & $\lambda_{I}$ & \multicolumn{1}{l@{}}{$\kappa_{I}$}\\
\hline
0000 & \emptyset&0.0444& 1.0000& 0\\
0001 & 4&0.0307& 0.5800& 0\\
0010 & 3&0.0307& 0.5800& 0\\
0011 & 34&0.0403& 0.3700& 0.0336\\
0100 & 2&0.0346& 0.6200& 0\\
0101 & 24&0.0323& 0.3724& 0.0128\\
0110 & 23&0.0323& 0.3724& 0.0128\\
0111 & 234&0.0547& 0.2422& -0.0020\\
1000 & 1&0.0482& 0.7000& 0\\
1001 & 14&0.0491& 0.4220& 0.0160\\
1010 & 13&0.0491& 0.4220& 0.0160\\
1011 & 134&0.0875& 0.2750& -0.0026\\
1100 & 12&0.0828& 0.4660& 0.0320\\
1101 & 124&0.0979& 0.2853& -0.0038\\
1110 & 123&0.0979& 0.2853& -0.0038\\
1111 & 1234&0.1875& 0.1875& 0.0006\\
\hline
\end{tabular}
\end{table}

If we have only tree cumulants $K\in\cM_{T}^{\kappa}$, we can still
identify the parameters of the model up to the label switching on the
inner nodes using Corollary~\ref{cor:formulas}.\vadjust{\goodbreak} Recall that if
$|I|\leq3,$ then $\kappa_{I}=\mu_{I}$ so, for example,
\begin{eqnarray*}
\bar{\mu}_{r}^{2}&=&\frac{\mu_{123}^{2}}{\mu_{123}^{2}+4\mu_{12}\mu
_{13}\mu_{23}}=0.36,
\\
\eta^{2}_{r,1}&=&\frac{\mu_{123}^{2}+4\mu_{12}\mu_{13}\mu_{23}}{\mu
_{23}^{2}}=0.25,
\\
\eta^{2}_{r,a}&=&\frac{\mu_{14}^{2}}{\mu_{12}^{2}}\frac{\mu
_{123}^{2}+4\mu_{12}\mu_{13}\mu_{23}}{\mu_{134}^{2}+4\mu_{13}\mu
_{14}\mu_{34}}=0.25.
\end{eqnarray*}
Note that the entries in Table \ref{tab:quartet1} can be computed in
several different ways. However, by Corollary~\ref{cor:formulas} this
does not matter. For instance, to compute $\bar{\mu}_{r}$ we picked
$1,2,3$ as three leaves separated by $r$. If, instead of $1,2,3,$ we
used $1,2,4,$ the answer would be the same since
\[
\bar{\mu}_{r}^{2}=\frac{\mu_{124}^{2}}{\mu_{124}^{2}+4\mu_{12}\mu
_{14}\mu_{24}}=0.36.
\]
Finally, in Appendix \ref{app:signs} we show that in this case we have
exactly four possible distinct choices for combinations of signs of
these parameters. The first one is the original one with all $\eta
_{u,v}>0,$ which we denote by $\omega$:
\begin{eqnarray*}
\eta_{r,1}&=&0.5, \qquad   \eta_{r,2}=0.4, \qquad   \eta_{r,a}=0.5, \qquad   \eta
_{a,3}=0.4, \qquad   \eta_{a,4}=0.4, \\
   \bar{\mu}_{r}&=&-0.6,  \qquad  \bar{\mu}_{a}=-0.4,\vadjust{\goodbreak}
\end{eqnarray*}
where we omit $\bar{\mu}_{i}$ for $i=1,2,3,4$ since these are
constant for all points in $\widehat{\Omega}_{T}$. We obtain three
remaining points by using local sign switching as defined in Appendix~%
\ref{app:signs}, which are $(\eta_{r,1},\eta_{r,2},\eta_{r,a},\eta
_{a,3},\eta_{a,4}, \bar{\mu}_{r}, \bar{\mu}_{a})=(-0.5,-0.4, -0.5,
0.4,0.4, 0.6, -0.4)$ or $(0.5,0.4, -0.5,\allowbreak  -0.4,-0.4, -0.6, 0.4)$
or $(-0.5,-0.4, 0.5, -0.4,-0.4, -0.6, -0.4)$.

\section{Discussion}

The reparametrization of Bayesian tree models with hidden variables
given herein has illuminated the structure of these tree models and has
enabled us to establish some identifiability results. However, the
applicability of the new coordinate system reaches far beyond
understanding identifiability. Some additional results will be
presented in forthcoming papers where we generalize both results of
\cite{auvray2006sad} and \cite{settimi1998gbg}, obtaining the full
semi-algebraic description of this model class, and results of \cite
{rusakov2006ams}, on the asymptotic approximation of the marginal
likelihood integrals.

The results given here can be extended in a straightforward way to the
case when all hidden variables are binary but all leaf variables are
arbitrary. It is less clear how the methods extend to tree models for
arbitrary finite discrete random variables, or more generally, to other
discrete graphical models. However, the extension to Gaussian models on
trees appears to be straightforward.

The definition of tree cumulants in (\ref{eq:kappa-in-rho}) can be
generalized using other posets than $\Pi_{T}$. This opens many
interesting possibilities to investigate more general coordinate
systems for binary models. They all share certain useful properties of
classical cumulants. In particular, Lemma \ref{lem:splitzero} is true
if the poset of tree partitions is replaced by any other lattice of
partitions. We will report on this result in a forthcoming paper.

\begin{appendix}
\renewcommand{\theequation}{\arabic{equation}}

\section{Change of coordinates}
\setcounter{equation}{30}
\subsection{From probabilities to central moments}\label{sec:change}

Let $\Delta_{2^{n}-1}$ be the set of all possible probability
distributions of a binary vector $X=(X_{1},\ldots, X_{n})$ as defined
in (\ref{eq:prob-simplex}). Let $\cC_{n}$ be the set of all possible
central moments $\mu_{I}$ for $I\in[n]_{\geq2}$ and means $\lambda
_{1},\ldots,\lambda_{n}$. In this section, we show that there exists
a polynomial isomorphism between $\Delta_{2^{n}-1}$ and $\cC_{n}$.

First, perform a change of coordinates from the raw probabilities
$p=[p_{\a}]$ to the non-central moments $\mathbf{\lambda}=[\lambda
_{\a}]$ for $\a=(\a_1,\ldots,\a_n)\in\{0,1\}^n$. This is a linear
map \mbox{$f_{p\lambda}\dvtx\R^{2^{n}}\rightarrow\R^{2^{n}}$}, where $\mathbf
{\lambda}=f_{p\lambda}(p)$ is defined as follows:
\begin{equation}\label{eq:lambda_in_p}
\lambda_\a=\sum_{\a\leq\beta\leq\mathbf{1}} p_\beta\qquad
\mbox{for any } \a\in\{0,1\}^{n},
\end{equation}
where $\mathbf{1}$ denotes the vector of ones and the sum is over all
binary vectors $\beta$ such that $\a\leq\beta\leq\mathbf{1}$ in
the sense that $\a_i\leq\beta_i\leq1$ for all $i=1,\ldots,n$. In
particular, $\lambda_{\mathbf{0}}=1$ for all probability\vadjust{\goodbreak}
distributions. Therefore, the image $\cL_{n}=f_{p\lambda}(\Delta
_{2^{n}-1})$ is contained in the hyperplane defined by $\lambda
_{\mathbf{0}}=1$. The map, $f_{p\lambda}\dvtx \Delta
_{2^{n}-1}\rightarrow\cL_{n}$, is invertible and hence we can obtain
coordinates on $\cL_{n}$ given by $\lambda_{\a}$ for all $\a\in\{
0,1\}^{n}$ such that $\a\neq\mathbf{0}$. The inverse of~$f_{p\lambda
}$ is the map, $f_{\lambda p}=f_{p\lambda}^{-1}\dvtx\cL_{n}\rightarrow
\Delta_{2^{n}-1}$, and is given by
\begin{equation}\label{eq:p-in-lambda}
p_\a=\sum_{\a\leq\beta\leq\mathbf{1}}(-1)^{|\beta-\a|} \lambda
_\beta\qquad\mbox{for $\a=(\a_1,\ldots,\a_n)\in\{0,1\}^n$}.
\end{equation}

The linearity of the expectation implies that the central moments can
be expressed in terms of non-central moments. In particular,
\begin{equation}\label{eq:central}
\mu_\a=\sum_{\mathbf{0}\leq\beta\leq\a}(-1)^{|\beta|} \lambda
_{\a-\beta} \prod_{i=1}^n \lambda_{e_i}^{\beta_i} \qquad\mbox{for }\a\in\{0,1\}^n,
\end{equation}
where $|\beta|=\sum_{i}\beta_{i}$. Using these equations, we can
transform variables from the non-central moments $[\lambda_{\a}]$ to
another set of variables given by all the means $\lambda
_{e_{1}},\ldots, \lambda_{e_{n}}$, where $e_1,\ldots, e_n$ are
standard basis vectors in $\R^n$, and central moments $[\mu_{\alpha
}]$ for $\a\in\{0,1\}^{n}$. The polynomial mapping $f_{\lambda\mu
}\dvtx\R^{2^{n}}\rightarrow\R^{{n}}\times\R^{2^{n}}$ is the identity
on the first $n$ variables corresponding to the means $\lambda
_{e_{1}},\ldots, \lambda_{e_{n}}$ and is defined by (\ref
{eq:central}) on the remaining variables. The image of $f_{\lambda\mu
}$ is contained in the subspace $\cH\subset\R^{n}\times\R^{2^{n}}$
given by $\mu_{e_{1}}=\cdots=\mu_{e_{n}}=0$.
It is easy to show (see, e.g., equation~(5), \cite
{balakrishnan1998nrb}) that the inverse of $f_{\lambda\mu}\dvtx\R
^{2^{n}}\rightarrow\cH$ is given as $f_{\mu\lambda}=f_{\lambda\mu
}^{-1}:\cH\rightarrow\R^{2^{n}}$ defined by
\begin{equation}\label{eq:noncentral}
\lambda_\a=\sum_{\mathbf{0}\leq\beta\leq\a}\mu_{\a-\beta}
\prod_{i=1}^n \lambda_{e_i}^{\beta_i} \qquad\mbox{for }\a\in\{
0,1\}^n.
\end{equation}

Let $\cC_{n}$ denote $f_{\lambda\mu} (\cL_{n})$. Then $\cC_{n}$ is
contained in $\cH$ and $\mu_{\mathbf{0}}=1$. We have, therefore,
obtained coordinates of $\cC_{n}$ given by $\lambda_{e_{1}},\ldots,
\lambda_{e_{n}}$ together with $\mu_{\a}$ for all $\a\in\{0,1\}
^{n}$ such that $|\a|\geq2$.

\vspace*{-3pt}\subsection{A reparametrization for general Markov
models}\vspace*{-3pt}
\label{sec:repar-models}
Let $T=(V,E)$ be a rooted tree with $n$ leaves and root $r$. Note that
for a tree $1+2|E|=|V|+|E|$ so the number of free parameters in (\ref
{eq:p_albar}) and (\ref{eq:p_albar2}) is $|V|+|E|$. We define a~polynomial map $f_{\theta\omega}:\R^{|V|+|E|}\rightarrow\R
^{|V|+|E|}$ from the original set of parameters of $\Theta_{T}$ given
by the root distribution and the conditional probabilities for each of
the edges to a set of parameters given as follows:
\begin{eqnarray}\label{eq:uij1}
\eta_{u,v}&=&\theta^{(v)}_{1|1}-\theta^{(v)}_{1|0} \qquad\mbox{for
each $(u,v)\in E$}  \quad \mbox{and}\nonumber
\\[-8pt]
\\[-8pt]
 \bar{\mu}_v&=&1-2\lambda_v\qquad\mbox{for each } v\in V,
\nonumber
\end{eqnarray}
where $\lambda_{v}=\E Y_{v}$ is a polynomial in the original
parameters $\theta$ of degree depending on the path from the root to
$v$. Let $(r,v_{1},\ldots, v_{k},v)$ be a directed path in $T$. Then
\[
\lambda_{v}=\sum_{\a\in\{0,1\}^{k+1}} \theta^{(v)}_{1|\a
_{k}}\theta^{(v_{k})}_{\a_{k}|\a_{k-1}}\cdots\theta^{(r)}_{\a_{r}}.\vadjust{\goodbreak}
\]

Let $\Omega_{T}=f_{\theta\omega}(\Theta_{T})$. The inverse map
$f_{\omega\theta}\dvtx \Omega_{T}\rightarrow\Theta_{T}$ has the
following form. For each edge $(u,v)\in E$ we have
\begin{eqnarray}\label{eq:theta-inverse}
\theta^{(v)}_{1|0}&=&\frac{1-\bar{\mu}_{v}}{2}-\eta_{u,v}\frac
{1-\bar{\mu}_{u}}{2},\nonumber
\\[-8pt]
\\[-8pt]
\theta^{(v)}_{1|1}&=&\frac{1-\bar{\mu}_{v}}{2}+\eta_{u,v}\frac
{1+\bar{\mu}_{u}}{2}
\nonumber
\end{eqnarray}
and $\theta^{(r)}_{1}=\frac{1-\bar{\mu}_{r}}{2}$.


\subsection{The non-degenerate case}\label{app:nondeg}

In this section, we derive the submodel of $\cM_{T}^{\kappa}=\psi
_{T}(\Omega_{T}),$ defined as the image of $\psi_{T}$ constrained to
the subset $\Omega_{T}^{0}$ of $\Omega_{T}$ given by $\bar{\mu
}_{v}^{2}<1$ for all $v\in V$. We define a smooth transformation on
$\Omega_{T}^{0}$ that enables us to change coordinates from $((\bar
{\mu}_{v}),(\eta_{u,v}))$ to $((\bar{\rho}_{v}),(\rho_{uv}))$, where
\begin{equation}\label{eq:zdef-rhouv}
\bar{\rho}_{v}=\frac{2\bar{\mu}_{v}}{\sqrt{1-\bar{\mu
}_{v}^{2}}},\qquad\rho_{uv}=\sqrt{\frac{1-\bar{\mu
}_{u}^{2}}{1-\bar{\mu}_{v}^{2}}}\eta_{u,v}.
\end{equation}
It is easily checked that this map is invertible since
\begin{equation}\label{eq:repar-to-rhos}
\bar{\mu}_{v}=\frac{\bar{\rho}_{v}}{\sqrt{4+\bar{\rho
}_{v}^{2}}},\qquad\eta_{u,v}=\sqrt{\frac{4+\bar{\rho
}_{u}^{2}}{4+\bar{\rho}_{v}^{2}}}\rho_{uv}.
\end{equation}
The inequality constraints defining $\Omega_{T}^{0}$ are given by
(\ref{eq:constraints}) and the fact that $\bar{\mu}_{v}\in(-1,1)$
for all $v\in V$. To express this in terms of the new coordinates, let
$t_{v}$ be defined by
\begin{equation}\label{eq:def-tefal}
t_{v}=\sqrt{1+ \biggl(\frac{\bar{\rho}_{v}}{2} \biggr)^{2}}+\frac
{\bar{\rho}_{v}}{2}\in(0,\infty).
\end{equation}
Then (\ref{eq:constraints}) becomes
\begin{eqnarray}\label{eq:new-constraints}
-t_{u}t_{v}&\leq&\rho_{uv}\leq\frac{t_{u}}{t_{v}},\nonumber
\\[-8pt]
\\[-8pt]
-\frac{1}{t_{u}t_{v}}&\leq&\rho_{uv}\leq\frac{t_{v}}{t_{u}}.
\nonumber
\end{eqnarray}

Transform the tree cumulants to a new coordinate system given by $\bar
{\rho}_{1},\ldots, \bar{\rho}_{n}$ and
\begin{equation}\label{eq:repar-to-rhos2}
\rho_{I}=\frac{2^{|I|}\kappa_{I}}{\prod_{i\in I}\sqrt{1-\bar{\mu
}_{i}^{2}}}\qquad\mbox{for all } I\in[n]_{\geq2},
\end{equation}
so that $\rho_{ij}$ is the correlation between $X_{i}$ and $X_{j}$.
The change of coordinates on $\Omega_{T}^{0}$ and~$\cK_{T}$ induces a
new parametrization of $\cM_{T}^{0}$. The parametrization\vadjust{\goodbreak} is given by
the identity on the first $n$ coordinates corresponding to $\bar{\rho
}_{i}$ for $i=1,\ldots,n$ and
\begin{equation}\label{eq:param-in-rhos}
\rho_{I}=\prod_{v\in V(I)\setminus I}\bar{\rho}_{v}^{\deg
(v)-2}\prod_{(u,v)\in E(I)}\rho_{uv}\qquad\mbox{for all } I\in
[n]_{\geq2}.
\end{equation}
In particular, each $\rho_{I}$ has an attractive monomial form. To
prove (\ref{eq:param-in-rhos}), simply substitute~(\ref
{eq:repar-to-rhos}) and (\ref{eq:repar-to-rhos2}) into (\ref
{eq:kappa_def_general}) to obtain
\[
\rho_{I} \prod_{i\in I}\frac{1}{\sqrt{4+\rho_{i}^{2}}}=\frac
{1}{4+\rho_{r(I)}^{2}}\prod_{v\in V(I)\setminus I} \biggl(\frac{\bar
{\rho}_{v}}{\sqrt{4+\bar{\rho}_{v}^{2}}} \biggr)^{\deg v-2}
 \prod_{(u,v)\in E(I)} \sqrt{\frac{4+\bar{\rho}_{u}^{2}}{4+\bar
{\rho}_{v}^{2}}}\rho_{uv}
\]
or, equivalently,
\begin{eqnarray*}
\rho_{I} &=&\prod_{v\in V(I)\setminus I}{\bar{\rho}_{v}}^{\deg
v-2}\prod_{(u,v)\in E(I)} \rho_{uv} \\
&&{}\times \frac{1}{4+\rho_{r(I)}^{2}} \prod_{v\in V(I)} \biggl(\frac
{1}{\sqrt{4+\bar{\rho}_{v}^{2}}} \biggr)^{\deg v-2}     \prod
_{v\in V(I)} \sqrt{\frac{4+\bar{\rho}_{\pa(v)}^{2}}{4+\bar{\rho
}_{v}^{2}}}.
\end{eqnarray*}
Next, we show that the term in the second line of the equation above is
equal to one. This follows from the fact that every $v\in V(I)$ apart
from the root is a parent of exactly $\deg(v)-1$ nodes and has one
parent; the root has no parents and is a parent of $\deg(r(I))$ nodes.

\section{\texorpdfstring{Proof of Proposition \protect\ref{prop:monomial}}
{Proof of Proposition 4.1}}\label{app:prop}
\setcounter{equation}{42}

It suffices to prove (\ref{eq:kappa_def_general}) for $I=[n]$ because
the general result for $I\subset[n]$ obviously follows by restriction
to the subtree $T(I)$ since each inner node of $T(I)$ has degree at
most three. The proof proceeds by induction with respect to the number
of leaves of $T$. First, we show that the result is true for $n=2$.
Since by definition $\kappa_{12}=\mu_{12}$ we need to prove
that\looseness=1
\begin{equation}\label{eq:prod-etas}
\mu_{12}=\frac{1}{4}({1-\bar{\mu}_{r}^{2}})\prod_{(u,v)\in E}\eta_{u,v},
\end{equation}\looseness=0
where $r$ is the root of $T$. If any of the nodes of $V$ represents a
degenerate random variable, then the global Markov properties in (\ref
{eq:GMP}) imply that $X_{1}\indep X_{2}$. In this case, the left-hand
side of (\ref{eq:prod-etas}) is zero. However, as we show next, one of
the factors on the right-hand side of~(\ref{eq:prod-etas}) must vanish
as well. We prove this by contradiction. Suppose that both $\bar{\mu
}_{r}^{2}\neq1$ and $\eta_{u,v}\neq0$ for all $(u,v)\in E$. By
Remark \ref{rem:phyl}, this implies that all the nodes of $T$
represent non-degenerate random variables, which leads to contradiction.

So assume now that every random variable in the system is
non-degenerate. From (\ref{eq:gen_rhoxy|h}), by taking $I=\{1\}$, $J=\{
2\},$ we have
\[
\mu_{12}=\tfrac{1}{4}({1-\bar{\mu}_{r}^{2}}) \eta_{r,1}\eta_{r,2}
\]
so it suffices to show that
\begin{eqnarray}\label{eq:pom-two-covar}
({1-\bar{\mu}_{r}^{2}})\eta_{r,1}&=&({1-\bar{\mu}_{r}^{2}})\prod
_{(u,v)\in E({r1})} \eta_{u,v}\quad\mbox{and}\nonumber
\\[-10pt]
\\[-10pt]
({1-\bar{\mu}_{r}^{2}})\eta_{r,2}&=&({1-\bar{\mu}_{r}^{2}})\prod
_{(u,v)\in E({r2})} \eta_{u,v}.
\nonumber\vspace*{-2pt}
\end{eqnarray}
If $r=1$ or $r$ is a parent of $1,$ then the first equation in (\ref
{eq:pom-two-covar}) is trivially satisfied. Assume that the length of
the path between $r$ and $1$ is greater than one. Let
$(r,h_{m},h_{m-1},\ldots,h_{1},1)$ be the directed path $E(r1)$
joining $r$ with $1$. Then, because $Y_{r}\indep Y_{1}|Y_{h_{1}}$, by
(\ref{eq:gen_rhoxy|h}) we have that
\begin{equation}\label{eq:IDEpom}
\tfrac{1}{4}({1-\bar{\mu}_{r}^{2}})\eta_{r,1}={\mu_{r1}}=\tfrac
{1}{4}({1-\bar{\mu}_{h_{1}}^{2}}) \eta_{h_{1},r}\eta_{h_{1},1}.\vspace*{-2pt}
\end{equation}
Similarly, because $Y_{r}\indep Y_{h_{k}}|Y_{h_{k+1}}$ for each
$k=1,\ldots,m-1,$ then again by (\ref{eq:gen_rhoxy|h})
\[
\tfrac{1}{4}(1-\bar{\mu}^{2}_{h_{k}})\eta_{h_{k},r}=\tfrac
{1}{4}(1-\bar{\mu}^{2}_{h_{k+1}})\eta_{h_{k+1},r}\eta_{h_{k+1},h_{k}}.\vspace*{-2pt}
\]
Substituting this expression for all subsequent $k=1,\ldots, m-1$ into
(\ref{eq:IDEpom}) we can now conclude that
\begin{equation}\label{eq:IDEpom2}
\tfrac{1}{4}({1-\bar{\mu}_{r}^{2}})\eta_{r,1}=\tfrac{1}{4}({1-\bar
{\mu}_{h_{m}}^{2}}) \eta_{h_{m},r}\eta_{h_{m},h_{m-1}}\cdots\eta
_{h_{2},h_{1}}\eta_{h_{1},1}.\vspace*{-2pt}
\end{equation}
But since $\frac{1}{4}({1-\bar{\mu}_{h_{m}}^{2}}) \eta
_{h_{m},r}=\mu_{r h_{m}}=\frac{1}{4}({1-\bar{\mu}_{r}^{2}}) \eta
_{r,h_{m}}$, equation (\ref{eq:IDEpom2}) implies that
\begin{equation}\label{eq:aux-product2}
({1-\bar{\mu}_{r}^{2}})\eta_{r,1}=(1-\bar{\mu}_{r}^{2})\prod
_{(u,v)\in E(r1)}\eta_{u,v}.\vspace*{-2pt}
\end{equation}
The second equation in (\ref{eq:pom-two-covar}) is proved simply by
changing the index from $1$ to $2$ above.


Now assume the proposition is true for all $k\leq n-1$ and let $T$ be a
tree with $n$ leaves. If one of the inner nodes of $T$ is degenerate,
then by the global Markov properties in (\ref{eq:GMP}) there exists an
edge split $C_{1}|C_{2}$ of the set of leaves such that
$X_{C_{1}}\indep X_{C_{2}}$. The left-hand side is zero by Lemma~\ref
{lem:splitzero}. Again, by Remark \ref{rem:phyl}, if both $\bar{\mu
}_{r}^{2}\neq1$ and $\eta_{u,v}\neq0$ for all $(u,v)\in E$, then
$\bar{\mu}_{v}^{2}\neq1$ for all $v\in V$. Hence, on the right-hand side
of equation (\ref{eq:prod-etas}), either $\bar{\mu}_{r}^{2}=1$ or
one of the $\eta_{u,{v}}$ vanishes. Consequently, (\ref
{eq:prod-etas}) is satisfied.

We assume now that all the inner nodes of $T$ represent non-degenerate
random variables. As $n\geq3$, we can always find two leaves separated
from all the other leaves by an inner node. We shall call such a pair
an \textit{extended cherry}. Denote the leaves by $1,2$ and the inner
node by $a$. Let $A=\{3,\ldots,n\}$ and let $T(aA)$ be the minimal
subtree of $T$ spanned $a\cup A$. Note that the global Markov
properties in (\ref{eq:GMP}) give that, for each $C\subseteq A,$ we
have $(X_1,X_2)\indep X_C|H_a$. Using (\ref{eq:gen_rhoxy|h}), we can
conclude that
\begin{equation}\label{eq:pom1}
\mu_{12C}=\mu_{12}\mu_{C}+\tfrac{1}{4}(1-\bar{\mu}_{a}^{2})\eta
_{a,12}\eta_{a,C}=\mu_{12}\mu_{C}+\eta_{a,12} \mu_{aC}.\vspace*{-2pt}
\end{equation}

Let $e\in E$ be the edge incident with $a$ separating $1$ and $2$ from
all other leaves, that is, such that $e$ induces the split $\nu
=12|\hat{1}_{A}$. For each $\pi\in\Pi_{T}$, if $\pi$ is induced by
removing $E_{\pi}\subset E,$ then ${\pi}\wedge\nu$ is induced
by
removing $E_{\pi}\cup e$. Let $\rho=12|\hat{0}_A\in\Pi_T$. Since
$\{1,2\}$ forms an extended cherry and all the inner nodes of $T$ have
degree at most three, it follows that $a$ necessarily has degree three
in $T$ and is a leaf of $T(aA)$. The\vadjust{\goodbreak} \textit{trimming map} with
respect to $\{1,2\}$ is the map $[\rho,\hat{1}]\rightarrow\Pi
_{T(aA)}$ such that $\pi\mapsto\widetilde{\pi}$ is defined by
changing the block $12C$ in $\pi\in[\rho,\hat{1}]$ to $aC$. Note
that the trimming map constitutes an isomorphism of posets between
$[\rho,\hat{1}]$ and $\Pi_{T(aA)}$.

It follows from the definition of tree cumulants in (\ref
{eq:kappa-in-rho}) that
\begin{equation}\label{eq:pom-kappa-split}
\kappa_{1\cdots n}=\sum_{\pi\in[\rho,\hat{1}]} \mm(\pi,\hat
{1})\prod_{B\in\pi}\mu_B+\sum_{\pi\notin[\rho,\hat{1}]} \mm
(\pi,\hat{1})\prod_{B\in\pi}\mu_B.\vspace*{-2pt}
\end{equation}
The second summand in (\ref{eq:pom-kappa-split}) is zero since every
$\pi\in\Pi_{T}$ such that $\pi\notin[\rho,\hat{1}]$ necessarily
contains either $1$ or $2$ as one of the blocks and $\mu_{1}=\mu_{2}=0$.
Applying (\ref{eq:pom1}) to each $\mu_{12C}$ for each $\pi\in[\rho
,\hat{1}],$ we obtain
\[
\prod_{B\in\pi}\mu_{B}=\prod_{B\in{\pi}\wedge\nu} \mu
_{B}+\eta_{a,12} \prod_{B\in\widetilde{\pi}}\mu_{B}\vspace*{-2pt}
\]
and hence
\begin{equation}\label{eq:mon_pom1}
\kappa_{1\cdots n}=\sum_{\pi\in[\rho,\hat{1}]} \mm(\pi,\hat
{1})\prod_{B\in{\pi}\wedge\nu}\mu_B+\eta_{a,12} \sum_{{\pi
}\in[\rho,\hat{1}]} \mm(\pi,\hat{1})\prod_{B\in\widetilde{\pi
}}\mu_B.\vspace*{-2pt}
\end{equation}

The first summand in (\ref{eq:mon_pom1}) can be rewritten as
\begin{equation}\label{eq:podsummand}
\sum_{\delta\in[\rho,\nu]} \biggl[  \biggl(\sum_{{\pi}\wedge\nu
=\delta}\mm(\pi,\hat{1}) \biggr)\prod_{B\in\delta}\mu_B
\biggr].\vspace*{-2pt}
\end{equation}
However, from Lemma \ref{lem:stanley}, since $\nu\neq\hat{1}$, for
each $\delta$ the sum $\sum_{{\pi}\wedge\nu=\delta}\mm(\pi,\hat
{1})$ in (\ref{eq:podsummand}) is zero. It follows that
\[
\kappa_{1\cdots n}=\eta_{a,12} \sum_{{\pi}\in[\rho,\hat{1}]}
\mm(\pi,\hat{1})\prod_{B\in\widetilde{\pi}}\mu_B.\vspace*{-2pt}
\]

By Proposition 4 in \cite{rota1964fct}, the M\"{o}bius function of
$[\rho,\hat{1}]$ is equal to the restriction of the M\"{o}bius
function on $\Pi_T$ to the interval $[\rho,\hat{1}]$. The trimming
map constitutes an isomorphism between $[\rho,\hat{1}]$ and $\Pi
_T(aA)$. Consequently, the M\"{o}bius function on $[\rho,\hat{1}]$ is
equal to the M\"{o}bius function on $\Pi_{T(aA)}$. It follows that
\begin{eqnarray*}
 \kappa_{1\cdots n}&=&\eta_{a,12}  \biggl(\sum_{{\pi}\in
[\rho,\hat{1}]} \mm(\pi,\hat{1})\prod_{B\in\widetilde{\pi}}\mu
_B \biggr)\\[-2pt]
 &=& \eta_{a,12} \biggl (\sum_{{\pi}\in\Pi_{T(aA)}} \mm
_{aA}(\pi,\hat{1}_{aA})\prod_{B\in{\pi}}\mu_B \biggr)=\eta
_{a,12} \kappa_{aA}.\vspace*{-2pt}
\end{eqnarray*}
Since $X_{1}\indep X_{2}|H_{a}$, by the second equation in Proposition
\ref{prop:indep}, $\eta_{a,12}=\bar{\mu}_{a}\eta_{a,1}\eta
_{a,2}$. Since $|aA|=n-1$, by using the induction assumption
\[
\kappa_{aA}=\frac{1}{4} \bigl(1-\bar{\mu}_{r(aA)}^{2} \bigr) \prod
_{v\in V(aA)\setminus aA} \bar{\mu}_{v}^{\deg(v)-2}\prod_{(u,v)\in
E(aA)} \eta_{u,v},\vspace*{-2pt}\vadjust{\goodbreak}
\]
where the degree is taken in $T(aA)$. We have two possible scenarios:
either $r(aA)\neq a$ or $r(aA)= a$. In the first case, $r(a1)=r(a2)=a$
and by (\ref{eq:aux-product2})
\[
\eta_{a,1}\eta_{a,2}=\prod_{(u,v)\in E(12)}\eta_{u,v}
\]
and hence
\begin{equation}\label{eq:kappa1-n}
\kappa_{1\cdots n}= \biggl(\bar{\mu}_{a}\prod_{(u,v)\in E(12)}\eta
_{u,v} \biggr)\kappa_{aA}.
\end{equation}
In the second case, either $r(a1)=a$ and $r(a2)=r$ or $r(a1)=r$ and
$r(a2)=a$ and so
\[
\eta_{a,1}\eta_{a,2}=\frac{1-\bar{\mu}_{r}^{2}}{1-\bar{\mu
}_{a}^{2}}\bar{\mu}_{a}\prod_{(u,v)\in E(12)}\eta_{u,v}.
\]
Hence,
\begin{equation}\label{eq:kappa1-n2}
\kappa_{1\cdots n}= \biggl(\frac{1-\bar{\mu}_{r}^{2}}{1-\bar{\mu
}_{a}^{2}}\prod_{(u,v)\in E(12)}\eta_{u,v} \biggr)\kappa_{aA}.
\end{equation}
The degree of $a$ in $T$ is three and the degree of all the other inner
nodes of $T(12)$ is two. Moreover, $E=E(aA)\cup E(12)$ and $V\setminus
[n]=(V(aA)\setminus aA)\cup(V(12)\setminus\{1,2\})$. It follows that
both (\ref{eq:kappa1-n}) and (\ref{eq:kappa1-n2}) satisfy (\ref
{eq:kappa_def_general}).

\section{Proofs of the theorems}\label{app:proofs}
\setcounter{equation}{53}

\begin{pf*}{Proof of Theorem \ref{lem:fibers}}
If each inner node of $T$ has degree at least three in $\widehat{T}$,
then for each inner node $u$ it is possible to find $i,j,k\in[n]$
separated by $u$ in $\widehat{T}$. So $\hat{\mu}_{ij}\hat{\mu
}_{ik}\hat{\mu}_{jk}\neq0$. Thus, by (\ref{eq:means}), we can
determine all values $\bar{\mu}_{u}^{2}=\hat{\mu}_{u}^{2}\neq1$.
Since, by Remark \ref{rem:steel_moulton}(ii), all the equivalence
classes in $[E\setminus\widehat{E}]$ are just single edges, we can
identify all $\eta_{u,v}^{2}=\hat{\eta}_{u,v}^{2}\neq0$ for all
$(u,v)\in E\setminus\widehat{E}$ by Lemma \ref{lem:steel_moulton}.

We now show that, because all equivalence classes in $[\widehat{E}]$
are singletons, $\eta_{w,w'}=0$ for every $(w,w')\in\widehat{E}$. By
construction, for each $(w,w')\in\widehat{E}$, either both $w$ and
$w'$ have degrees at least three in $\widehat{T}$ or one of them is a
leaf and the other has degree at least three in $\widehat{T}$.
Therefore, there exist $i,j\in[n]$ such that $E(ij)\cap\widehat{E}=\{
(w,w')\}$ by the construction of $\widehat{E}$. We have that $\hat
{\mu}_{ij}=0$. However, $\eta_{u,v}=\hat{\eta}_{u,v}\neq0$ for all
$(u,v)\in E\setminus\widehat{E}$. Because $\bar{\mu
}_{r(ij)}^{2}=\hat{\mu}_{r(ij)}^{2}\neq1,$ it follows by (\ref
{eq:muij}) that $\eta_{w,w'}=0$. Therefore, the values of all the
parameters are fixed up to signs and in this case $\widehat{\Omega
}_{T}$ is finite. The proof that there are exactly $2^{|V|-n}$ points
in this fiber is provided in Appendix~\ref{app:signs}.

To prove the second statement of Theorem \ref{lem:fibers}, first note
that, since every inner node of $T$ has degree at least two in
$\widehat{T}$, it follows by Lemma \ref{lem:non-deg} that for each
$v\in V$, \mbox{$\bar{\mu}_{v}^{2}<1$}. This implies that the\vadjust{\goodbreak} $\hat
{p}$-fiber lies in $\Omega_{T}^{0}\subset\Omega_{T}$ as defined in
Appendix \ref{app:nondeg}. We can apply a smooth transformation over
this subset to a second space $\Omega'_{T}\subseteq\R^{|V|+|E|}$
whose coordinates are given by $\bar{\rho}_{v}$ for $v\in V$ and
$\rho_{uv}$ for $(u,v)\in E$. The map is defined by (\ref
{eq:zdef-rhouv}) and is invertible with the inverse defined in (\ref
{eq:repar-to-rhos}).

To investigate the geometry of the $\hat{p}$-fiber in $\Omega'_{T}$,
first list all the defining constraints. For all $i=1,\ldots,n$ we
have that $\bar{\mu}_{i}=\hat{\mu}_{i}$ because $\hat{p}$
determines the sample means of the observed nodes. Hence the value of
$\bar{\rho}_{i}$ is determined as well. Write $\bar{\rho}_{i}=\hat
{\rho}_{i}$ for all $i=1,\ldots,n$, where $\hat{\rho}_{i}$ is the
image of $\hat{\mu}_{i}$ under (\ref{eq:zdef-rhouv}). For each inner
node $v$ whose degree in $\widehat{T}$ is at least three, we can find
$i,j,k\in[n]$ separated in $\widehat{T}$ by $v$. The value of $\bar
{\mu}_{v}^{2}$ is determined by (\ref{eq:means}), which is well
defined because $\hat{\mu}_{ij}\hat{\mu}_{ik}\hat{\mu}_{jk}>0$.
Therefore, the value of~$\bar{\rho}_{v}^{2}$, for each~$v$ whose
degree in $\widehat{T}$ is at least three, is fixed $\bar{\rho
}_{v}^{2}=\hat{\rho}_{v}^{2}$, where $\hat{\rho}_{v}^{2}=\frac
{4\hat{\mu}_{v}^{2}}{1-\hat{\mu}_{v}^{2}}$
by~(\ref{eq:zdef-rhouv}).\looseness=-1

Next, we show that for every $(u,v)\in\widehat{E}$ we must have that
$\rho_{uv}=0$. This follows by essentially the same argument as in the
first part of the proof. Because the degrees of both $u$ and $v$ are at
least two, there exist $i,j\in[n]$ such that $E(ij)\cap\widehat{E}=\{
(u,v)\}$. In particular, $\hat{\mu}_{ij}=0$ and so by (\ref
{eq:muij}) $\eta_{u,v}=0$. Moreover, for any path $E(kl)$ in
$[E\setminus\widehat{E}]$ the value of $\rho_{kl}^{2}$ is constant
by Lemma~\ref{lem:steel_moulton}. So write $\rho_{kl}=\hat{\rho
}_{kl}$. By (\ref{eq:param-in-rhos}), we have that
\begin{equation}\label{eq:rho-kl-hat}
\hat{\rho}_{kl}=\prod_{(u,v)\in E(kl)}\rho_{uv}.
\end{equation}
Finally, for any degree-two node $v$ the parameter $\bar{\rho}_{v}$
can take any real value and each $\rho_{uv}$ is constrained to satisfy
(\ref{eq:new-constraints}). This completes the list of constraints
defining the image of the $\hat{p}$-fiber in $\Omega'_{T}$.

We now show that this image is diffeomorphic to a union of polyhedra.
Let $\rho=((\bar{\rho}_{v}),(\rho_{uv}))$ be any point in the
transformed $\hat{p}$-fiber. Then $\rho$ lies in a linear subspace~$\cL$ of $\R^{|V|+|E|}$ given by $\rho_{uv}=0$ for all $(u,v)\in
\widehat{E}$. Since $\rho_{uv}\neq0$ for all $(u,v)\in E\setminus
\widehat{E}$, we can define the following further smooth change of
coordinates on $\cL$. Let $s\dvtx E\rightarrow\{-1,0,1\}$ be any possible
sign assignment for $(\rho_{uv})$ such that $s(u,v)=\operatorname{sgn}(\rho
_{uv})$ and $\operatorname{sgn}(\rho_{ij})=\prod_{(u,v)\in E(ij)}s(u,v)$ for
all $i,j\in[n]$ (cf. Appendix \ref{app:signs}). Then $s$ induces an
open orthant~$\R^{|E\setminus\widehat{E}|}_{s}$ defined by
$s(u,v)\rho_{uv}> 0$ for all $(u,v)\in E\setminus\widehat{E}$.
Moreover, the disjoint union of $\cU_{s}=\R^{|V|}\times\R
_{s}^{|E\setminus\widehat{E}|}\subset\cL$, for all possible sign
assignments $s$, covers the $\hat{p}$-fiber, that is, each point of
the $\hat{p}$-fiber lies in one of the $\cU_{s}$. Note also that on
each $\cU_{s}$ the sign of~$\bar{\rho}_{v}$ for all nodes of the
degree at least three is fixed. This follows from the fact that
by~(\ref{eq:param-in-rhos})\looseness=-1
\[
\rho_{ijk}=\bar{\rho}_{v}\prod_{(u,w)\in E(ijk)}\rho_{uw},
\]\looseness=0
for any three leaves $i,j,k\in[n]$ separated by $v$ in $\widehat{T}$.
Since on each $\cU_{s}$ the signs of $\rho_{uw}$ for all $(u,w)\in
E(ijk)$ are fixed, the sign of $\bar{\rho}_{v}$ also has to be fixed
to match the sign of~$\rho_{ijk}$. We write $\bar{\rho}_{v}=\hat
{\rho}_{v}^{s}$ on $\cU_{s}$.

On each $\cU_{s}$ define a map to the space $\R^{|V|+|E\setminus
\widehat{E}|}$ with coordinates given by $\nu_{uv}$ for $(u,v)\in
E\setminus\widehat{E}$ and $z_{v}$ for $v\in V$. The map is a
diffeomorphism defined as follows. We set
\[
\nu_{uv}=\log(s(u,v)\rho_{uv})\qquad\mbox{for all }(u,v)\in
E\setminus\widehat{E}.\vadjust{\goodbreak}
\]
Next, for every $v\in V$ we substitute $\bar{\rho}_{v}$ for $t_{v}$
as defined in (\ref{eq:def-tefal}). This is an invertible
transformation because
\[
\bar{\rho}_{v}=\frac{t_{v}^{2}-1}{t_{v}},
\]
which is well defined since $t_{v}>0$ for all $v\in V$. We then simply
substitute $t_{v}$ for $z_{v}=\log t_{v}$.

In this new coordinate system, the $\hat{p}$-fiber restricted to $\cU
_{s}$ is a union of polyhedra. The defining constraints are as follows. First,
\begin{eqnarray}\label{eq:str-polyh}
z_{i}&=& \hat{z}_{i} \qquad  \mbox{for all leaves }
i=1,\ldots,n,\nonumber
\\[-10pt]
\\[-10pt]
z_{v}&=& \hat{z}_{v}^{s}  \qquad  \mbox{for all $v$ with degree at least three
in }\widehat{T}.
\nonumber
\end{eqnarray}
Here, $\hat{z}_{i},\hat{z}_{v}^{s}$ are real numbers obtained as
images of $\hat{\rho}_{i}$, $\hat{\rho}_{v}^{s}$, respectively.
Moreover, for each $E(kl)\in[E\setminus\widehat{E}]$
\begin{equation}\label{eq:nuuv-sum}
\sum_{(u,v)\in E(kl)}\nu_{uv}=\log|\hat{\rho}_{kl}|
\end{equation}
subject to additional inequality constraints
\begin{eqnarray}\label{eq:ineqonompr}
\nu_{uv}&\leq&\min\{z_{u}-z_{v},z_{v}-z_{u}\} \qquad  \mbox{if } s(u,v)=1,\nonumber\\[-2pt]
\nu_{uv}&\leq&\min\{z_{u}+z_{v},-z_{u}-z_{v}\} \qquad  \mbox{if } s(u,v)=-1,
\mbox{ for each $(u,v)\in E\setminus\widehat{E}$ and}\\[-2pt]
z_{v}&>&0  \qquad  \mbox{for the inner nodes of degree $2$}.\nonumber
\end{eqnarray}
These inequalities follow from (\ref{eq:new-constraints}). Since all
these constraints are linear, they define a~polyhedron in $\R
^{|V|+|E\setminus\widehat{E}|}$. Therefore the $\hat{p}$-fiber is a
disjoint union of subsets each of which is diffeomorphic to a polyhedron.

To show the dimension of each polyhedron is equal to $2l_{2}$, we must
ensure that the dimension of the smallest affine subspace containing
this polyhedron is $2l_{2}$. Since $z_{v}>0$ for all $v\in V$ it is
easily checked that the inequalities in (\ref{eq:ineqonompr}) do not
induce any equality. Therefore, the description of the affine span is
obtained from the description of the polyhedron (given by (\ref
{eq:str-polyh})--(\ref{eq:ineqonompr})) by
suppressing all inequalities in (\ref{eq:ineqonompr}). The dimension
of the ambient space is $|V|+|E\setminus\widehat{E}|$; the
codimension is given by the number of equations in (\ref
{eq:str-polyh}) and (\ref{eq:nuuv-sum}). Hence the codimension is
equal to $|V|-l_{2}+|[E\setminus\widehat{E}]|$. For each $E(kl)\in
[E\setminus\widehat{E}]$ one has that $|E(kl)|-1$ is equal to the
number of degree-two nodes in $E(kl)$. By summing over all $E(kl)$ it
follows that $|E\setminus\widehat{E}|-|[E\setminus\widehat
{E}]|=l_{2}$. Therefore, the dimension of the polyhedron is given by
\[
(|V|+|E\setminus\widehat{E}|)-(|V|-l_{2}+|[E\setminus\widehat{E}]|)=2l_{2}.
\]
Since the dimension of the affine span of a polyhedron is equal to its
dimension, the dimension is equal to $2l_{2}$ as required.
\end{pf*}

\begin{pf*}{Proof of Theorem \ref{prop:singular}}
Let $V_{0}\subseteq\widehat{V}$ and $E_{0}\subseteq\widehat{E}$ and
\begin{equation}\label{eq:omega-AB}
\Omega_{(V_0,E_0)} =\{\omega\in{\Omega}_{T}\dvt   \bar{\mu
}_{v}^{2}=1 \mbox{ for all } v\in V_{0},   \eta_{u,v}=0\mbox{ for
all } (u,v)\in E_{0} \}.\vadjust{\goodbreak}
\end{equation}
We say that $(V_0,E_0)$ is \textit{minimal for} $\widehat{\Sigma}$
if for every point $\omega$ in $\Omega_{(V_0,E_0)}$ and for every~$i,j\in[n]$ such that $\hat{\mu}_{ij}=0$ we have that $\mu
_{ij}(\omega)=0$ and furthermore that $(V_0,E_0)$ is minimal with such
a property (with respect to inclusion on both coordinates).

To illustrate the motivation behind this definition, consider the
tripod tree singular case in Example \ref{ex:sing}. If $T$ is rooted
in the inner node, we have four minimal subsets of $2^{\widehat
{V}}\times2^{\widehat{E}}$: $(\{h\},\emptyset)$, $(\emptyset,\{
(h,1),(h,2)\})$, $(\emptyset,\{(h,1),(h,3)\})$ and $(\emptyset,\{
(h,2),(h,3)\})$.

We now show that the $\hat{p}$-fiber satisfies
\begin{equation}\label{eq:p-fiber-decomp}
\widehat{\Omega}_{T}=\bigcup_{(V_0,E_0)\min\!.} \Omega
_{(V_0,E_0)}\cap\widehat{\Omega}_{T}.
\end{equation}
The first inclusion ``$\subseteq$'' follows from the fact that if
$\omega\in\widehat{\Omega}_{T}$, then $\mu_{ij}(\omega)=\hat{\mu
}_{ij}$ for all~$i,\allowbreak j\in[n]$. In particular, $\mu_{ij}(\omega)=0$
whenever $\hat{\mu}_{ij}=0$. Therefore, $\omega\in\Omega
_{(V_0,E_0)}\cap\widehat{\Omega}_{T}$ for $(V_0,E_0)$ minimal. The
second inclusion is obvious.

%
%


For each minimal $(V_0,E_0)$ the set $\Omega_{(V_0,E_0)}\cap\widehat
{\Omega}_{T}$ is a union of disjoint manifolds in~$\R^{|V|+|E|}$
constrained to $\Omega_{T}$. To show this, consider first all the
connected components $T_{i}=(V_{i},E_{i})$ for $i=1,\ldots, k$ of
$\widehat{T}$ except isolated inner nodes of $\widehat{T}$. By Remark~%
\ref{rem:steel_moulton}(iv), all these components are trees with a
set of leaves contained in $[n]$. The projection of the parameter space
$\Omega_{T}$ to the parameters for the marginal model $\cM
_{T_{i}}^{\kappa}$ is denoted by~$\Omega_{i}$. It is therefore a
projection of $\Omega_{T}$ on $\bar{\mu}_{v}$ for $v\in V_{i}$ and
$\eta_{u,v}$ for $(u,v)\in E_{i}$. By Theorem~\ref{lem:fibers}, each
component $T_{i}$ induces a manifold with corners in $\Omega_{i}$,
denoted by $\widehat{\Omega}_{i}$. Hence there exists a manifold
$M_{i}$ in $\R^{|V_{i}|+|E_{i}|}$ such that $\widehat{\Omega
}_{i}=M_{i}\cap\Omega_{i}$. The constraints on the remaining
coordinates are given by: $\bar{\mu}_{v}^{2}=1$ for all $v\in V_{0}$
and $\eta_{u,v}=0$ for $(u,v)\in E_{0}$. These algebraic equations
define a union $M_{(V_0,E_0)}$ of affine subspaces in $\R^{|\widehat
{V}|+|\widehat{E}|}$ with coordinates given by $\bar{\mu}_{v}$ for
$v\in\widehat{V}$ and $\eta_{u,v}$ for $(u,v)\in\widehat{E}$.

For each $(V_{0},E_{0}),$ consider the union of manifolds $M\subset\R
^{|V|+|E|}$ given as the Cartesian product of $M_{(V_0,E_0)}$ and
$M_{i}$ for $i=1,\ldots, k$. The restriction of $M$ to $\Omega_{T}$
is exactly $\Omega_{(V_0,E_0)}\cap\widehat{\Omega}_{T}$. Now we
have that
\begin{equation}\label{eq:int-AB}
\bigcap_{(V_0,E_0) \min\!.}  \bigl(M_{(V_0,E_0)}\times
M_{1}\times\cdots\times M_{k} \bigr)= \biggl(\bigcap_{(V_0,E_0)
\min\!.} M_{(V_0,E_0)} \biggr)\times M_{1}\times\cdots\times M_{k}.
\end{equation}
However, $\bigcap_{(V_0,E_0) \min\!.} M_{(V_0,E_0)}$ is equal to
\[
\bigl\{\omega\in\R^{|V|+|E|}\dvt  \bar{\mu}_{v}^{2}=1\mbox{ for all }
v\in\widehat{V}, \eta_{u,v}=0 \mbox{ for all } (u,v)\in\widehat
{E}\bigr\},
\]
where, after the restriction to $\Omega_{T}$, the intersection in
(\ref{eq:int-AB}) is equal to the deepest singularity.
\end{pf*}

\section{Sign patterns for parameters}\label{app:signs}

Let $\hat{p}\in\cM_{T}$ such that each inner node of $T$ has degree
at least three in the corresponding forest $\widehat{T}$. By the proof
of Theorem \ref{lem:fibers}, there is a finite number of points
$\theta\in\Theta_{T}$ such that $f_{T}(\theta)=\hat{p}$. By
definition, this set of points is denoted by $\widehat{\Theta}_{T}$.
Corollary \ref{cor:formulas} gives the formulae for the parameters
modulo signs, which suggests that $|\widehat{\Theta
}_{T}|=2^{|V|+|E|}$. However, not all sign choices are possible. Let
$m$ be the number of inner nodes of $T$. We will show that the number
of possible choices of signs is, in fact, equal to $2^{m}$, that is,
$|\widehat{\Theta}_{T}|=2^{m}$. We also show how to obtain all the
points in $\widehat{\Theta}_{T}$ given one of them. This construction
becomes especially simple when expressed in the new parameters defined
by~(\ref{eq:uij1}).\looseness=1

Let $\theta$ be a point in $\widehat{\Theta}_{T}$ ($\widehat{\Theta
}_{T}$ is finite and non-empty) and let $\omega=f_{\theta\omega
}(\theta)$. We assign signs to each edge of $T$ using the map $s\dvtx
E\rightarrow\{-1,0,1\}$ such that for every
$(u,v)\in E$, $s(u,v)=\operatorname{sgn}(\eta_{u,v})$, where $\eta_{u,v}$ are parameters in $\omega$.
Let $h$ be an inner node of $T$. On $\widehat{\Omega}_{T}$ we define
the operation of local sign switching $\delta_{h}$ such that $\delta
_{h}(\omega)= \omega'$ where $\eta'_{u,v}=-\eta_{u,v}$ if one of
the ends of $(u,v)$ is in $h$ and $\eta'_{u,v}=\eta_{u,v}$ otherwise;
$\bar{\mu}_{h}'=-\bar{\mu}_{h}$ and $\bar{\mu}_{v}'=\bar{\mu
}_{v}$ for all $v\neq h$. We have that $\bar{\mu}_{i}'=\bar{\mu
}_{i}$ and hence $\lambda_{i}'=\lambda_{i}$ for all leaves
$i=1,\ldots,n$. Let now $I\in[n]_{\geq2}$. Then, from (\ref
{eq:kappa_def_general}),
\[
\kappa_{I}(\omega')=\frac{1}{4}\bigl(1-\bar{\mu}_{r(I)}^{2}\bigr)\prod
_{v\in V(I)\setminus I} (\bar{\mu}_{v}')^{\deg(v)-2} \prod
_{(u,v)\in E(I)}\eta_{u,v}'.
\]
We have two cases: either $h$ lies in $V(I)$ or not. In the first case,
\[
\kappa_{I}(\omega')=(-1)^{\deg(h)-2}(-1)^{\deg(h)}\kappa
_{I}(\omega)=\kappa_{I}(\omega).
\]
In the second case, $\omega'=\omega$ and hence trivially $\kappa
_{I}(\omega')=\kappa_{I}(\omega)$. It follows that $\omega'\in
\widehat{\Omega}_{T}$ and therefore the operator $\delta
_{h}\dvtx\widehat{\Omega}_{T}\rightarrow\widehat{\Omega}_{T}$ is well
defined. The local sign switchings form a~group $\cG$ that is
isomorphic to the multiplicative group $Z_{2}^{m}$. By composing
distinct local switchings we obtain $2^{m}$ different points in
$\widehat{\Omega}_{T}$. Hence the orbit of $\omega$ in $\widehat
{\Omega}_{T}$ has exactly~$2^{m}$ elements.

It remains to show that there are no other orbits of $\cG$ in
$\widehat{\Omega}_{T}$. Let $\omega\in\widehat{\Omega}_{T}$ and
let $\omega'$ be a point in $\Omega_{T}$ such that $({\eta
_{u,v}')}^{2}={\eta_{u,v}}^{2}$ for all $(u,v)\in E$ and ${(\bar{\mu
}_{v}')}^{2}={\bar{\mu}_{v}}^{2}$ for all inner nodes $v$ of $T$,
which is a necessary condition for $\omega'$ to be in $\widehat
{\Omega}_{T}$. Assume that $\omega'$ is not in the orbit of $\omega
$. We will show below that this implies that $\omega'$ cannot lie in
the $\hat{p}$-fiber. It will then follow that the orbit of $\omega$
constitutes the whole $\widehat{\Omega}_{T}$ and hence $|\widehat
{\Omega}_{T}|=2^{m}$.

We proceed by contradiction. Thus, let $\omega'\in\widehat{\Omega
}_{T}$ and we want to show that $\omega'=\delta(\omega)$ for some
$\delta\in\cG$. Since $\omega$ can be replaced by any other point
in its orbit, we can assume that $\operatorname{sgn}(\bar{\mu}_{v})=
\operatorname{sgn}(\bar{\mu}'_{v})$ for all $v\in V$. Since $\omega,\omega'\in
\widehat{\Omega}_{T}$, for every $i,j,k\in[n]$ by (\ref{eq:kappa_def_general})
applied for~$\kappa_{ij}$ and $\kappa_{ijk}$,
respectively, we have that
\[
\prod_{(u,v)\in E(ij)} s(u,v)=\prod_{(u,v)\in E(ij)} s'(u,v),\qquad
\prod_{(u,v)\in E(ijk)} s(u,v)=\prod_{(u,v)\in E(ijk)} s'(u,v).
\]
It follows that $\prod_{(u,v)\in E(vi)}s(u,v)=\prod_{(u,v)\in
E(vi)}s'(u,v)$ for each inner node $v$ and leaf $i$. It immediately
implies that $s(u,v)=s'(u,v)$ for all $(u,v)\in E$ and hence $\omega
=\omega'$. In this way we have shown that $\omega'$ is in the orbit
of $\omega$ under $\cG$.
\end{appendix}

\section*{Acknowledgement}
We are very grateful to a referee whose
extensive comments enabled us to substantially improve this
paper.


\printhistory


\begin{thebibliography}{21}

\bibitem{allman2009identifiability}
\begin{barticle}[mr]
\bauthor{\bsnm{Allman},~\bfnm{Elizabeth~S.}\binits{E.S.}},
  \bauthor{\bsnm{Matias},~\bfnm{Catherine}\binits{C.}} \AND
  \bauthor{\bsnm{Rhodes},~\bfnm{John~A.}\binits{J.A.}}
(\byear{2009}).
\btitle{Identifiability of parameters in latent structure models with many
  observed variables}.
\bjournal{Ann. Statist.}
\bvolume{37}
\bpages{3099--3132}.
\bid{doi={10.1214/09-AOS689}, issn={0090-5364}, mr={2549554}}
\end{barticle}
\endbibitem

\bibitem{auvray2006sad}
\begin{bmisc}[auto:STB|2011-03-03|12:04:44]
\bauthor{\bsnm{Auvray},~\bfnm{Vincent}\binits{V.}},
  \bauthor{\bsnm{Geurts},~\bfnm{Pierre}\binits{P.}} \AND
  \bauthor{\bsnm{Wehenkel},~\bfnm{Louis}\binits{L.}}
(\byear{2006}).
\bhowpublished{A semi-algebraic description of discrete naive Bayes models with two
  hidden classes.
In \textit{Proc. Ninth International Symposium on Artificial Intelligence
  and Mathematics, Fort Lauderdale, Florida}. Available at \url{http://anytime.cs.umass.edu/aimath06/}.}
\end{bmisc}
\endbibitem

\bibitem{balakrishnan1998nrb}
\begin{barticle}[mr]
\bauthor{\bsnm{Balakrishnan},~\bfnm{N.}\binits{N.}},
  \bauthor{\bsnm{Johnson},~\bfnm{Norman~L.}\binits{N.L.}} \AND
  \bauthor{\bsnm{Kotz},~\bfnm{Samuel}\binits{S.}}
(\byear{1998}).
\btitle{A note on relationships between moments, central moments and cumulants
  from multivariate distributions}.
\bjournal{Statist. Probab. Lett.}
\bvolume{39}
\bpages{49--54}.
\bid{doi={10.1016/S0167-7152(98)00027-3}, issn={0167-7152}, mr={1649335}}
\end{barticle}
\endbibitem

\bibitem{chang1996frm}
\begin{barticle}[mr]
\bauthor{\bsnm{Chang},~\bfnm{Joseph~T.}\binits{J.T.}}
(\byear{1996}).
\btitle{Full reconstruction of {M}arkov models on evolutionary trees:
  Identifiability and consistency}.
\bjournal{Math. Biosci.}
\bvolume{137}
\bpages{51--73}.
\bid{doi={10.1016/S0025-5564(96)00075-2}, issn={0025-5564}, mr={1410044}}
\end{barticle}
\endbibitem

\bibitem{feller1971ipt}
\begin{bbook}[mr]
\bauthor{\bsnm{Feller},~\bfnm{William}\binits{W.}}
(\byear{1971}).
\btitle{An Introduction to Probability Theory and Its Applications. {V}ol.
  {II}},
\bedition{2nd} ed.
\baddress{New York}: \bpublisher{Wiley}.
\bid{mr={0270403}}
\end{bbook}
\endbibitem

\bibitem{geiger2001sef}
\begin{barticle}[mr]
\bauthor{\bsnm{Geiger},~\bfnm{Dan}\binits{D.}},
  \bauthor{\bsnm{Heckerman},~\bfnm{David}\binits{D.}},
  \bauthor{\bsnm{King},~\bfnm{Henry}\binits{H.}} \AND
  \bauthor{\bsnm{Meek},~\bfnm{Christopher}\binits{C.}}
(\byear{2001}).
\btitle{Stratified exponential families: Graphical models and model selection}.
\bjournal{Ann. Statist.}
\bvolume{29}
\bpages{505--529}.
\bid{doi={10.1214/aos/1009210550}, issn={0090-5364}, mr={1863967}}
\end{barticle}
\endbibitem

\bibitem{lauritzen96}
\begin{bbook}[mr]
\bauthor{\bsnm{Lauritzen},~\bfnm{Steffen~L.}\binits{S.L.}}
(\byear{1996}).
\btitle{Graphical Models}.
\bseries{Oxford Statistical Science Series}
\bvolume{17}.
\baddress{Oxford}: \bpublisher{Clarendon Press}.
\bid{mr={1419991}}
\end{bbook}
\endbibitem

\bibitem{mccullagh1987tms}
\begin{bbook}[mr]
\bauthor{\bsnm{McCullagh},~\bfnm{Peter}\binits{P.}}
(\byear{1987}).
\btitle{Tensor Methods in Statistics}.
\baddress{London}: \bpublisher{Chapman \& Hall}.
\bid{mr={0907286}}
\end{bbook}
\endbibitem

\bibitem{moulton2004ppo}
\begin{barticle}[mr]
\bauthor{\bsnm{Moulton},~\bfnm{Vincent}\binits{V.}} \AND
  \bauthor{\bsnm{Steel},~\bfnm{Mike}\binits{M.}}
(\byear{2004}).
\btitle{Peeling phylogenetic `oranges'}.
\bjournal{Adv. in Appl. Math.}
\bvolume{33}
\bpages{710--727}.
\bid{doi={10.1016/j.aam.2004.03.003}, issn={0196-8858}, mr={2095862}}
\end{barticle}
\endbibitem

\bibitem{pearltarsi86}
\begin{barticle}[mr]
\bauthor{\bsnm{Pearl},~\bfnm{Judea}\binits{J.}} \AND
  \bauthor{\bsnm{Tarsi},~\bfnm{Michael}\binits{M.}}
(\byear{1986}).
\btitle{Structuring causal trees}.
\bjournal{J. Complexity}
\bvolume{2}
\bpages{60--77}.
\bid{issn={0885-064X}, mr={0925434}}
\end{barticle}
\endbibitem

\bibitem{rota1964fct}
\begin{barticle}[mr]
\bauthor{\bsnm{Rota},~\bfnm{Gian-Carlo}\binits{G.C.}}
(\byear{1964}).
\btitle{On the foundations of combinatorial theory. {I}. {T}heory of {M}\"obius
  functions}.
\bjournal{Probab. Theory Related Fields}
\bvolume{2}
\bpages{340--368}.
\bid{mr={0174487}}
\end{barticle}
\endbibitem

\bibitem{rotacumulants}
\begin{barticle}[mr]
\bauthor{\bsnm{Rota},~\bfnm{Gian-Carlo}\binits{G.C.}} \AND
  \bauthor{\bsnm{Shen},~\bfnm{Jianhong}\binits{J.}}
(\byear{2000}).
\btitle{On the combinatorics of cumulants}.
\bjournal{J. Combin. Theory Ser. A}
\bvolume{91}
\bpages{283--304}.
\bid{doi={10.1006/jcta.1999.3017}, issn={0097-3165}, mr={1779783}}
\end{barticle}
\endbibitem

\bibitem{rusakov2006ams}
\begin{barticle}[mr]
\bauthor{\bsnm{Rusakov},~\bfnm{Dmitry}\binits{D.}} \AND
  \bauthor{\bsnm{Geiger},~\bfnm{Dan}\binits{D.}}
(\byear{2005}).
\btitle{Asymptotic model selection for naive {B}ayesian networks}.
\bjournal{J. Mach. Learn. Res.}
\bvolume{6}
\bpages{1--35 (electronic)}.
\bid{issn={1532-4435}, mr={2249813}}
\end{barticle}
\endbibitem

\bibitem{semple2003pol}
\begin{bbook}[mr]
\bauthor{\bsnm{Semple},~\bfnm{Charles}\binits{C.}} \AND
  \bauthor{\bsnm{Steel},~\bfnm{Mike}\binits{M.}}
(\byear{2003}).
\btitle{Phylogenetics}.
\bseries{Oxford Lecture Series in Mathematics and Its Applications}
\bvolume{24}.
\baddress{Oxford}: \bpublisher{Oxford Univ. Press}.
\bid{mr={2060009}}
\end{bbook}
\endbibitem

\bibitem{settimi1998gbg}
\begin{bincollection}[auto:STB|2011-03-03|12:04:44]
\bauthor{\bsnm{Settimi},~\bfnm{Raffaella}\binits{R.}} \AND
  \bauthor{\bsnm{Smith},~\bfnm{Jim~Q.}\binits{J.Q.}}
(\byear{1998}).
\btitle{On the geometry of Bayesian graphical models with hidden variables}.
In \bbooktitle{UAI}
(\beditor{\bfnm{Gregory~F.}\binits{G.F.}~\bsnm{Cooper}} \AND
  \beditor{\bfnm{Moral}\binits{M.}~\bsnm{Seraf{\'{\i}}n}}, eds.)
\bpages{472--479}.
\baddress{San Francisco}: \bpublisher{Morgan Kaufmann}.
\end{bincollection}
\endbibitem

\bibitem{settimi2000gma}
\begin{barticle}[mr]
\bauthor{\bsnm{Settimi},~\bfnm{Raffaella}\binits{R.}} \AND
  \bauthor{\bsnm{Smith},~\bfnm{Jim~Q.}\binits{J.Q.}}
(\byear{2000}).
\btitle{Geometry, moments and conditional independence trees with hidden
  variables}.
\bjournal{Ann. Statist.}
\bvolume{28}
\bpages{1179--1205}.
\bid{doi={10.1214/aos/1015956712}, issn={0090-5364}, mr={1811324}}
\end{barticle}
\endbibitem

\bibitem{speed1983cumulants}
\begin{barticle}[mr]
\bauthor{\bsnm{Speed},~\bfnm{T.~P.}\binits{T.P.}}
(\byear{1983}).
\btitle{Cumulants and partition lattices}.
\bjournal{Austral. J. Statist.}
\bvolume{25}
\bpages{378--388}.
\bid{issn={0004-9581}, mr={0725217}}
\end{barticle}
\endbibitem

\bibitem{speicher94}
\begin{barticle}[mr]
\bauthor{\bsnm{Speicher},~\bfnm{Roland}\binits{R.}}
(\byear{1994}).
\btitle{Multiplicative functions on the lattice of noncrossing partitions and
  free convolution}.
\bjournal{Math. Ann.}
\bvolume{298}
\bpages{611--628}.
\bid{doi={10.1007/BF01459754}, issn={0025-5831}, mr={1268597}}
\end{barticle}
\endbibitem

\bibitem{spiegelhalter1993bae}
\begin{barticle}[mr]
\bauthor{\bsnm{Spiegelhalter},~\bfnm{David~J.}\binits{D.J.}},
  \bauthor{\bsnm{Dawid},~\bfnm{A.~Philip}\binits{A.P.}},
  \bauthor{\bsnm{Lauritzen},~\bfnm{Steffen~L.}\binits{S.L.}} \AND
  \bauthor{\bsnm{Cowell},~\bfnm{Robert~G.}\binits{R.G.}}
(\byear{1993}).
\btitle{Bayesian analysis in expert systems}.
\bjournal{Statist. Sci.}
\bvolume{8}
\bpages{219--283}.
\bnote{With comments and a~rejoinder by the authors}.
\bid{issn={0883-4237}, mr={1243594}}
\end{barticle}
\endbibitem

\bibitem{stanley2006enumerative}
\begin{bbook}[auto:STB|2011-03-03|12:04:44]
\bauthor{\bsnm{Stanley},~\bfnm{Richard~P.}\binits{R.P.}}
(\byear{2002}).
\btitle{Enumerative Combinatorics. Volume I. Cambridge
  Studies in Advanced Mathematics}
  \bvolume{49}.
\baddress{Cambridge}: \bpublisher{Cambridge Univ. Press}.
\end{bbook}
\endbibitem

\end{thebibliography}
\end{document}